\documentclass[11pt]{article}

\usepackage[letterpaper,margin=1in]{geometry} 

\usepackage[english]{babel}

\usepackage[
  style=apa, 
  bibstyle=authoryear,
  autocite=inline,
  maxcitenames=2,
 maxbibnames=8,
  backend=biber
]{biblatex}
\usepackage{setspace}
\usepackage[utf8]{inputenc}
\usepackage[T1]{fontenc}
\usepackage[english]{babel}
\usepackage{csquotes}  
\usepackage{graphicx}
\usepackage[export]{adjustbox}
\graphicspath{ {./images/} }

\usepackage{amsmath,amssymb}
\usepackage{mathabx,mathrsfs}
\usepackage[version=4]{mhchem}
\usepackage{stmaryrd}
\usepackage{bbold}

\usepackage{indentfirst}
\usepackage{caption}
\usepackage{float}
\usepackage{hyperref}
\usepackage{algorithm}
\usepackage{algorithmic}
\usepackage{mathtools}
\usepackage{array}
\usepackage{blindtext}
\usepackage{titlesec}
\usepackage{xcolor}
\usepackage[final]{microtype} 
\definecolor{darkblue}{HTML}{0B2F6A}
\usepackage{amsthm} 

\newtheorem{theorem}{Theorem} 
\newtheorem{proposition}[theorem]{Proposition} 

\hypersetup{colorlinks,linkcolor={blue},citecolor={blue},urlcolor={red}}

\newtheorem{remark}{Remark}

\DeclareMathOperator*{\Bern}{Bernoulli}
\DeclareMathOperator*{\STD}{StDev}
\newcommand{\Dev}[1]{\mathrm{Dev}(#1)}
\DeclareMathOperator*{\arginf}{arg\ inf}

\newcommand{\Xmax}{X_{\text{max}}}
\DeclareRobustCommand{\QQ}{\texorpdfstring{$Q$}{Q}}
\DeclareRobustCommand{\BB}{\texorpdfstring{$B$}{B}}
\newcommand{\Imax}{I_{\max}}
\newcommand{\Imin}{I_{\min}}
\newcommand{\Bmin}{B_{\min}}
\newcommand{\Bmax}{B_{\max}}
\newcommand{\EE}{\mathbb{E}}

\newcommand{\PP}{\mathbb{P}}

\newcommand{\bx}{\mathbf{x}}

\newcommand{\bC}{\mathbf{C}}
\newcommand{\bI}{\mathbf{I}}

\newcommand{\beginsupplement}{%
    \setcounter{section}{0}
    \renewcommand{\thesection}{S\arabic{section}}
    \renewcommand{\thesubsection}{S\arabic{section}.\arabic{subsection}}
    \setcounter{figure}{0}
    \renewcommand{\thefigure}{S\arabic{figure}}
    \setcounter{table}{0}
    \renewcommand{\thetable}{S\arabic{table}}
    \setcounter{equation}{0}
    \renewcommand{\theequation}{S\arabic{equation}}
}

\title{\LARGE \bf
Intraday Battery Dispatch for \\ Hybrid Renewable Energy Assets
}

\author{%
Thiha Aung\thanks{Department of Statistics and Applied Probability, University of California, Santa Barbara. Email: taung@ucsb.edu} \quad
and \quad
Mike Ludkovski\thanks{Department of Statistics and Applied Probability, University of California, Santa Barbara. Email: ludkovski@pstat.ucsb.edu}}
\date{}
\begin{document}
\maketitle
\begin{abstract}
We develop a mathematical model for intraday dispatch of co-located wind-battery energy assets. Focusing on the primary objective of firming grid-side actual production vis-a-vis the  preset day-ahead hourly generation targets, we conduct a comprehensive study of the resulting stochastic control problem across different firming formulations and wind generation dynamics. Among others, we provide a closed-form solution in the special case of a quadratic objective and linear dynamics, as well as design a novel adaptation of a Gaussian Process-based Regression Monte Carlo algorithm for our setting. Extensions studied include an asymmetric loss function for peak shaving, capturing the cost of battery cycling, and the role of battery duration. In the applied portion of our work, we calibrate our model to a collection of 140+ wind-battery assets in a synthetic Texas grid, benchmarking the economic benefits of firming based on outputs of a realistic unit commitment and economic dispatch solver. 
\end{abstract}


\section{Introduction}

Unlike conventional thermal generators that are fully dispatchable, variable renewable energy causes increased operational costs for power grids. Over- and under-generation requires deployment of regulation down/up reserves and is often associated with financial losses for the asset owner as well. Hybrid assets---defined as a  ``generating resource that is comprised of multiple generation or energy storage technologies controlled as a single entity behind a single point of interconnection'' \parencite{nerc2021gridforming}---have been gaining popularity as a means to mitigate the variability inherent to standalone renewable generation \parencite{Ahlstrom_2021}. According to \textcite{LBNL} report, 46\% of new renewable projects in the US listed in grid interconnection queues at the end of 2023 were of hybrid type, and total deployed hybrid capacity includes 49GW of renewable generation and 24.2 GWh of energy capacity, which is almost half of all deployed battery capacity in the US. 

A hybrid resource combines a renewable energy source---most commonly solar or wind---with an energy storage system. We focus on fully hybrid resources that included both co-located and co-controlled resources that are treated as a single grid participant. 
Multiple storage technologies are being utilized (battery storage, mechanical gravity, compressed air, hydrogen fuel cells) with lithium-ion battery energy storage system (BESS) being the most common at present. Thanks to the controlled storage system, hybrid assets offer additional flexibility, closing the gap to traditional dispatchable generators. The comprehensive report by \parencite{ESIGrep} lists more than half a dozen different drivers for resource hybridization, including avoided distribution and transmission upgrades, reduced curtailment, hedge against varying market conditions, and simplified procurement for utility generation off-takes. In turn, \textcite{EIA-2024} reports that hybrid resources utilize their batteries for many diverse use cases, from energy price arbitrage and ancillary service provision, to system peak shaving. 

A critical use case of hybrid asset concerns renewable firming/curtailment mitigation. This involves dispatching the battery to maximally match the hybrid output with an exogenous dispatch target. For instance, this target can correspond to the day-ahead generation forecast employed by the system operator for security-constrained unit commitment. Firming objectives also arise naturally in the context of resource adequacy analysis and Power Purchase Agreements (PPAs) offtakes that reward predictability of generation. 

In this paper, we investigate the operation of hybrid assets on operational time scales, focusing on optimal dispatch of the coupled BESS. We consider wind-hybrid assets with our first task corresponding to the firming objective. The BESS is dynamically controlled to counter the deviations of realized wind power output relative to a given dispatch target. Taking into account the intertemporal constraints of the BESS storage capacity and power rating limits, we study the resulting state-constrained  finite-horizon stochastic optimal control problem.  

Our developments span methodological, algorithmic and empirical aspects. Methodologically, we use a stochastic differential equation to directly model renewable energy generation and formulate a generic intra-day control problem for hybrid BESS dispatch. Starting with asset firming, we demonstrate its applicability for peak shaving and congestion/curtailment avoidance tasks. Moreover, we show that we can incorporate battery degradation control criteria into our flexible framework. Algorithmically, we directly implement the dynamic programming equation with state constraints by utilizing a simulation-based machine learning algorithm of Regression Monte Carlo (RMC)-type. Using actor-critic framework, we build two Gaussian Process (GP) emulators for the continuation-value and optimal control maps. Our algorithm is agnostic to assumed state dynamics which can be highly nonlinear and also agnostic to the objective function. We leverage two main features of the RMC: a statistical surrogate to approximate the continuation-value function and a customized simulation design. We propose direct optimization of the control over the cost functional, which includes the continuation-value surrogate.  

To provide empirical validation,  we carry out a comprehensive case study, calibrating our model to realistic wind power production data from the synthetic Texas-7k transmission grid, a digital twin of the Electric Reliability Council of Texas (ERCOT) grid. After fitting our probabilistic model for each wind unit in Texas-7k, we run grid-scale dispatch simulations to assess the impact of hybridization on asset firming and daily dispatch savings. Our simulations rely on actual (rather than model-based) wind and load profiles, providing an out-of-sample  quantification of how well our dispatch algorithm performs in real-dispatch conditions.

The literature on stochastic control of hybrid assets are still nascent and can be grouped in terms of the assumed control strategy types and objectives.
\textcite{johnson_partial_2017} assume a Bang–Bang control strategy for BESS to optimize daily commitments, approximating constraints via boundary conditions and solving the HJB PDE numerically.
\textcite{tankov_wind_BESS} likewise model a hybrid wind asset for arbitrage, using Bang–Bang control to track dispatch targets while introducing a decision variable for intraday energy trading.
In contrast, we optimize BESS control for dispatch tracking beyond Bang–Bang strategies, handling constraints without approximations through GP emulators. Neural networks have also seen as emulators to solve stochastic control problems through  deep BSDE  \parencite{han2020deep}, PDE  \parencite{Sirignano_2018}, and actor-critic  \parencite{hure2020applications} techniques. Alternatively, see \textcite{Carvalho_Energy_Trading} and \textcite{Focker_Planck_wind_BESS} for numerical methods for solving stochastic control of BESS without machine learning emulators. Other control frameworks applied to BESS in hybrid assets include Model Predictive Control \parencite{deterministic_Wind_BESS_MPC, Stochastic_MPC_Wind_BESS}, scenario-based stochastic programming \parencite{ZAREOSKOUEI20151105,DING2012571}, and Approximate Dynamic Programming \parencite{Belloni_2016}.
 

In addition to hybrid configurations, BESSs are also utilized for standalone applications, including energy arbitrage in real-time and/or day-ahead markets \parencite{zheng2022energy}, and participation in ancillary service markets \parencite{Arteaga_2019,
Dheepak_2018}. While BESSs operate in standalone/hybrid configurations, degradation is a key component of dispatch optimization. In optimal control, degradation function is quadratically approximated, with additional state variables introduced to capture/control cycling depth \parencite{MPCDegration, LEE2022107795}. In contrast, the standard rainflow method \parencite{multi-factormodel} is path-dependent and non-Markovian, making it unsuitable for standard DP formulations.

To tackle our stochastic control problem, we build on the Stochastic Hybrid Asset Dispatch Optimization with Gaussian Processes (SHADOw-GP) algorithm developed in \textcite{Aung_Ludkovski}. Compared to neural networks, GPs offer superior sample efficiency, requiring far fewer training points for accurate approximations; see \textcite{ludkovski2018simulation,alasseur2018regression,BALATA2021640} for related applications in renewable energy control. While \textcite{Aung_Ludkovski} focused on algorithmic development using simplified linear dynamics and quadratic objectives in stylized settings, this work advances the framework toward realistic modeling and deployment. Specifically, we adopt time-dependent Jacobi-type SDEs with non-Gaussian noise for wind dynamics, incorporate nonlinear objectives such as bang-bang control, curtailment mitigation, and degradation costs, and derive a novel analytical solution for optimal BESS control under “soft” State-of-Charge and control constraints. Furthermore, we highlight the practical impact of hybrid assets through an empirical case study using the synthetic Texas-7k grid.

The rest of the paper is organized as follows. Section \ref{sec:2} presents the stochastic optimal control problem set-up for the hybrid asset. A linear quadratic approximation of the constraints and analytical solution is discussed in Section \ref{sec:3}. In Section \ref{sec:4}, we present our SHADOw-GP algorithm. Section \ref{sec:5} is devoted to an extended case study of hybridizing wind assets in the synthetic Texas-7k grid. Results on additional optimization objectives are provided in Section \ref{sec:6}.
Section \ref{sec:8} outlines future research directions.

\section{Problem Formulation}\label{sec:2}
We consider a standard probability space $\left(\Omega, \mathcal{F},\left(\mathcal{F}_t\right)_{t\in [0,T]}, \PP \right)$. We work with two Markovian state variables: the exogenous wind power generation process $(X_t)_{t\in[ 0,T]}$ and the controlled State of Charge (SoC) process $(I_t)_{t\in[ 0,T]}$, both adapted to the filtration $\mathbb{F}=\left(\mathcal{F}_t\right)_{t\in [0,T]}$. The SoC process is controlled by the charge/discharge rate process $(B_t)_{t\in[ 0,T]}$, a continuous $\mathbb{F}$-adapted decision variable. The BESS is dynamically controlled throughout the $24$-hour operation period as $(X_t)$ evolves, subject to control and state constraints of the SoC process. This is interpreted as real-time adjustment by the hybrid asset owner; note that we are not referencing any power prices or bidding/auctions.

\subsection{Wind Power Generation Dynamics}

Let $X = (X_t)_{t \in [0,T]}$ denote the wind power generation process in MW. We model $(X_t)$ using a mean-reverting SDE, where the mean-reversion level $(m_t)_{t \in [0,T]}$ represents the day-ahead forecast and is assumed to be non-negative, $m_t \geq 0$. Because generation is non-negative and capped by the nameplate generation capacity, the domain of $X$ is bounded, $X_t \in [0, X_{\max}]$. To this end, we work with Jacobi diffusion SDE
\parencite{SDE_I}:
\begin{equation}\label{eq:OU}
    dX_t=  \alpha_t ( m_t-X_t) \, dt+ \sigma_t \sqrt{X_t(\Xmax-X_t)}\, dW_t, \qquad X_0=x_0,
\end{equation} 
where 
$\sigma_t >0$ is the volatility parameter, 
$\alpha_t > 0$ is the dimensionless mean-reversion parameter, and $ (W_t)_{t\in[0,T]}$ is a Wiener process, capturing the real-time generation uncertainty. Note that $\alpha_t$ and $\sigma_t$ are time-dependent to capture the diurnal variability in wind power output. 
\begin{remark}
A similar framework applies to PV–BESS or PV–wind–BESS hybrids. Unlike wind, solar generation is limited to daylight hours, so the domain of nonzero $X_t$ is a proper subset of $[0,T]$. Moreover, PV–BESS systems primarily aim to shift generation from morning to evening, rather than perform firming as considered here.
\end{remark}
\subsection{State-of-Charge Dynamics}

Let $I=(I_t)_{t\in[0,T]}$ denote the SoC process of stored energy in units MWh. $(I_t)$ is Markovian, subject to charge/discharge rate constraints and capacity constraints, with dynamics:
\begin{equation}\label{eq:I}
    dI_t =  \bigl(\eta B_{t} \mathbb{1}_{ \{B_t  > 0\} } +\frac{1}{\eta} B_{t}\mathbb{1}_{\{B_t <0\} } \bigr) \, dt, \qquad I_0 = \iota,
\end{equation}
where $(B_t)_{t\in[ 0,T]}$ is the controlled charge/discharge rate process in units MW, a continuous decision variable. $B_t>0$ denotes charging, while $B_t<0$ represents discharging. $\eta \in (0,1)$ denotes the dimensionless charge/discharge efficiency parameter that captures charging inefficiencies, which means that the battery dissipates more energy than it discharges and accumulates less energy than it charges \parencite{Kordonis_BESS_control,Stochastic_MPC_PV_BESS}. Typical values for $\eta$ are around 90-95\%. The SoC bounds are given by the battery capacity constraints:
\begin{equation} \label{SoCcap}
     I_{\min} \equiv \text{SoC}_{\min} \cdot I_{\text{cap}} \leq I_{t} \leq \text{SoC}_{\max} \cdot I_{\text{cap}} \equiv I_{\max}
\end{equation}
where $I_{\text{cap}}$ is the rated battery capacity in units MWh and $0\leq \text{SoC}_{\min}<\text{SoC}_{\max}\leq 1$ are SoC limits. Typical values for the latter are $\text{SoC}_{\min} = 0.05, \text{SoC}_{\max} = 0.95$. The BESS power rating is characterized by $B_{\max}>0$ and $B_{\min}<0$ that restrict the range of the control $B_t$:
\begin{equation}
\label{Bconstriants}
\begin{split}
    B_{\min} &\leq B_t \leq B_{\max}. 
\end{split}
\end{equation}
The Markovian structure of $(X_t,I_t)$ allows us to express the optimal control in feedback form $B_t = {\Psi(t,X_t,I_t)}$, adapted to filtration $\mathbb{F}$.

\subsection{Stochastic Control Problem }\label{SOCP}

\begin{figure}[ht]
    \centering
    \includegraphics[width=0.5\linewidth,trim=8in 14in 0in 8in,clip=true]{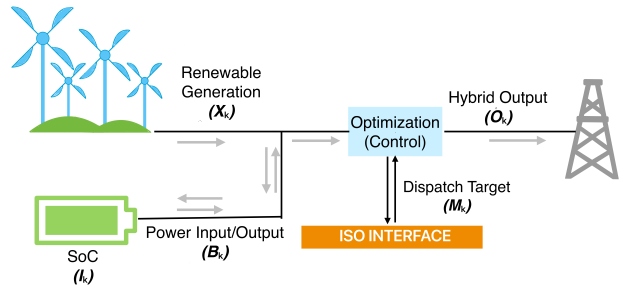} 
    
    \caption{A schematic description of optimizing the output $(O_k)$ of a hybrid wind asset.}
    \label{fig:example}
\end{figure}

A hybrid wind asset engages in grid dispatch by receiving dispatch signals $(M_t)_{t\in [0,T]}$ from the Independent System Operator (ISO). 
The  realized power output $X_t$ deviates from the dispatch target $M_t$ and the hybrid asset charges/discharges the battery so that the net output $O_t:= X_t-B_t$ is close to $ M_t$. That is, during operation the BESS agent dynamically decides the dispatch rate in MW as a function of the SoC and the wind power output; see Figure \ref{fig:example}. 

Mathematically, denote by $X^{t,x} = (X^{t,x}_s)_{s\in [t,T]}$ the solution of the SDE \eqref{eq:OU} starting from $x$ at $t$ and by $I^{B;t,\iota} = (I^{B;t,\iota}_s)_{s\in [t,T]}$ the solution of the SoC process \eqref{eq:I}, starting from $\iota$, and controlled by $B := (B_s)_{s\in [t,T]}$. Let $\mathcal{B}(t,x,\iota)$ be the set of feasible controls for the given initial condition:
\begin{equation}\label{eq:feasible_set}
    \mathcal{B}(t,x,\iota) = \left\{ B : B_s \text{ is } \mathcal{F}_s\text{-adapted}, \, B_{\text{min}} \leq B_s \leq B_{\text{max}}, \, \Imin \leq I_s \leq \Imax ,\ \forall s \in [t,T]\right\}.
\end{equation}
The value function $V:[t, T] \times \mathcal{D} \rightarrow \mathbb{R}$ of our control problem is :
\begin{equation}\label{eq:HJB}
V(t, x,\iota):=\inf _{B\in {\mathcal{B}(t,x,i)} } \EE \left[\int_t^T f(X^{t,x}_s,B_s,M_s) \,d s+g\left(I^{B;t,\iota}_T\right)\right],
\end{equation}
where the domain is the bounded set $\mathcal{D}:=\mathcal{X} \times \mathcal{I} \subset \mathbb{R}^2$ 
with $\mathcal{X} = [0,X_{\max}]$ and $\mathcal{I} = [I_{\min},I_{\max}]$. 
The running cost $f(\cdot)$ represents the stepwise performance criterion for optimal BESS dispatch over the interval $[0,T]$ which depends on the target $M_t$. The motivating  stepwise cost is a quadratic firming criterion
    \begin{align}\label{eq:tildeF}
    \tilde{f}(X,B,M) := (X-B-M)^2.
    \end{align}
The terminal cost $g(\cdot)$ captures the final SoC level that the BESS needs to end up at to be ready for the next operational day. 

\begin{remark} The natural firming target would be $M_t = \EE[X_t]$ which is different from the mean-reversion level $m_t$ as the dynamics of $(X_t)$ are time-dependent. 
\end{remark}

\subsection{Explicit Solution via Linear-Quadratic Approximation}\label{sec:3}

To set a baseline for latter developments, we first consider a linearized setup that admits an explicit solution. To this end, we cast \eqref{eq:HJB} into the linear-quadratic framework and then analytically solve the corresponding HJB equation. 

Dynamic programming \parencite{Pham_2009} suggests that the value function of \eqref{eq:HJB} can be characterized as the viscosity solution of the HJB PDE 
\begin{equation}\label{constrained_HJB}
\begin{aligned}
\partial_t V + \inf_{b \in \mathcal{B}_t} \Bigl\{ 
& \alpha_t (m_t - x)\, \partial_X V 
+ \bigl(\eta b \mathbb{1}_{\{b > 0\}} + \tfrac{1}{\eta} b \mathbb{1}_{\{b < 0\}}\bigr) \partial_I V \\
& + \tfrac{1}{2} \sigma_t^2 x (\Xmax - x) \partial_{XX} V 
+ f(x,b,M_t) 
\Bigr\} = 0.
\end{aligned}
\end{equation}
where $(x,{\iota}) \in \mathcal{D}$ with terminal condition $V(T,x,{\iota}) = g({\iota})$. 
Above, \begin{equation}\label{constraint_set_t}
    \mathcal{B}_t :=[ \Bmin \cdot \mathbb{1}_{I_t > I_{\min}}, \Bmax \cdot  \mathbb{1}_{I_t < I_{\max}} 
    ] 
\end{equation}
denotes the $I_t$-dependent admissible charge/discharge rates at time $t$. 
To solve our problem analytically, we relax the hard constraints in $\mathcal{B}_t$ using a  quadratic penalty on the state and control variables. Furthermore, we take a quadratic terminal cost and linearize the dynamics of $I_t$:
\begin{itemize}
    \item Linear dynamics for the SoC process by considering $\eta = 1$ ($100 \%$ efficient) BESS: 
\begin{equation}\label{eq:linearI}
    dI_t = B_t \, dt, \qquad I_0 = {\iota},
\end{equation}
    \item Quadratic stepwise cost $\bar{f}$, with the firming criterion modified to approximate the  hard state and control constraints with penalties  $c_1,c_2>0$. We let $I_m := \frac{I_{\text{max}}-I_{\text{min}}}{2}$ and take:
\begin{equation}\label{running cost}
     \bar{f}(X_t,B_t,M_t,I_t) = \tilde{f}(X_t,B_t,M_t)  + c_1 B_t^2+ c_2 (I_t-I_m)^2.
\end{equation}

    \item A quadratic terminal penalty centered around a target SoC level $I_{\text{target}}$:
\begin{equation}\label{eq:terminal_cost}
    \bar{g}(I_T) = \mathcal{P} \cdot (I_T-I_{\text{target}})^2 \qquad \text{with } \quad \mathcal{P} >0.
\end{equation}
\end{itemize}

\noindent Under assumptions \eqref{eq:OU}, \eqref{eq:linearI}-\eqref{eq:terminal_cost}, we obtain an HJB PDE with unconstrained real-valued control set:
\begin{multline}\label{unconstrained_HJB}
\partial_t \bar{V}+\inf_{\bar{b}\in \mathbb{R}} \Bigl\{\alpha_t  (m_t-x)\partial_X \bar{V} +\bar{b}\partial_{I}\bar{V}+\frac{1}{2} \sigma_t^2 x (\Xmax -x)  \partial_{XX} \bar{V} \\ +(x- \bar{b}-M_t)^2 + c_1 \bar{b}^2+ c_2 ({\iota}-I_m)^2\Bigr\}=0 
\end{multline}  with terminal condition $V(T,x,\iota) = \bar{g}(\iota)$. 

\begin{remark}

The underlying control problem spans a multi-day horizon in real-time grid dispatch.  However, solving a multi-day stochastic control problem is computationally expensive and furthermore, we do not know the dispatch targets beyond day-ahead. As a result, we instead opt for solving a sequence of one-day sub-problems. Since the true terminal condition (i.e., the value function at the start of day \( D+1 \)) is unknown at the end of day \( D \), we impose a \emph{heuristic terminal condition} given by \eqref{eq:terminal_cost}. This terminal cost serves as a \emph{soft guidance mechanism} to discourage excessive deviation of the terminal state from the targeted SoC  \( I_{\text{target}} \). Specifically, the terminal state at the end of day \( D \), denoted \( I_T \), becomes the initial state for the following day \( D+1 \). This ensures inter-day consistency in the absence of explicit multi-day coordination.
\end{remark}
\begin{proposition}\label{lq_control}
Assume that the solution \(\bar{V}\) of \eqref{unconstrained_HJB} is $\mathcal{C}^{1,2}([0,T] \times \mathcal{D})$. When the mean-reversion level is constant, $m_t = m$ 
 $\forall t\in [0,T]$, the value function is explicitly given by: 
\begin{equation}
\begin{aligned}
\bar{V}(t, x, \iota) &= P_1(t)(\iota - I_m)^2 + P_2(t)(\iota - I_m)(x - m) + P_3(t)(x - m)^2 \\
&\quad + P_4(t)(\iota - I_m) + P_5(t)(x - m) + P_6(t),
\end{aligned}
\end{equation}
where the functions $P_1, \ldots, P_6$ solve the system of backward Ricatti ordinary differential equations 
{\small
\begin{align} \left\{ 
\begin{aligned}
&\dot{P_1}(t) + c_2 - \kappa (P_1(t))^2 = 0, & P_1(T) = \mathcal{P}\\
&\dot{P_2}(t) -\alpha_t P_2(t) -\kappa P_1(t) P_2(t) + 2 \kappa P_1(t) = 0, & P_2(T) = 0 \\
&\dot{P_3}(t)+(1-\kappa) -(2 \alpha_t + \sigma^2_t) P_3(t) +\kappa P_2(t) -\frac{\kappa}{4} (P_2(t))^2 = 0, & P_3(T)=0 \\
&\dot{P_4}(t) + 2\kappa(m-M_t) P_1(t) - \kappa P_1(t) P_4(t) = 0, &  P_4(T) = 2\mathcal{P}(I_m - \iota_{\text{target}})\\
&\dot{P_5}(t) + 2(1-\kappa) (m-M_t)-\alpha_t P_5(t) + \kappa P_4(t) + \kappa(m-M_t) P_2(t)  \\
&\quad -\frac{\kappa}{2} P_2(t)P_4(t) + (\sigma^2_t \Xmax - 2 \sigma^2_t m) P_3(t)= 0, & P_5(T) = 0 \\
&\dot{P_6}(t) + (1-\kappa) (m-M_t)^2 + \kappa (m-M_t) P_4(t) -\frac{\kappa}{4} (P_4(t))^2 = 0, & P_6(T) = \mathcal{P}(I_m - \iota_{\text{target}})^2
\end{aligned} \right. \label{eq:ricatti}
\end{align}
}
with $\kappa = \frac{1}{1+c_1}$. 

Moreover, the optimal control in \eqref{unconstrained_HJB} is of the closed-loop feedback form  and is given by:
\begin{equation}\label{eq:b_bar_exact}
    \bar{b}(t,x,\iota) = \kappa (x-M_t) -\kappa P_1(t) (\iota-I_m) -\frac{\kappa}{2}P_2(t) (x-m) -\frac{\kappa}{2}P_4(t).
\end{equation}
\end{proposition}

\begin{proposition}\label{exists_n_unique}
Let $c_2, \kappa, \mathcal{P} > 0$. Assume $\alpha_t, \sigma_t, M_t \in L^1([0,T])$. Then the Riccati ODE system \eqref{eq:ricatti} admits a unique global solution on $[0,T]$ for any $T>0$.
\end{proposition}
Consequently, the existence and uniqueness of the optimal control in \eqref{eq:b_bar_exact} are guaranteed. Proofs of Propositions~\ref{lq_control} and~\ref{exists_n_unique} are provided in Section~S1 of the Supplementary Material.



\section{Numerical Algorithm }\label{sec:4}

We present our SHADOw-GP algorithm to solve the stochastic control problem (\ref{eq:HJB}. The algorithm is presented in discrete time with step size $\Delta t$ and $K$ total steps, such that $T=K \Delta t$. Denoting the time index at $t_k = k\Delta t$ with subscript $k$, the discrete-time versions of the wind power generation and SoC dynamics are

\begin{equation}\label{eq:OU_discrete}
        X_{k+1}=X_{k}+\alpha_k (m_k - X_k)\Delta t+\sigma_k  \sqrt{X_k (\Xmax-X_k)} \sqrt{\Delta t}\cdot Z_{k}, \qquad Z_k \sim \mathcal{N}(0, 1);
\end{equation}
\begin{equation}\label{eq:I_discrete}
    I_{k+1} = I_{k} + \bigl(\eta B_{k} \mathbb{1}_{ \{B_k  > 0\} } +\frac{1}{\eta} B_{k}\mathbb{1}_{\{B_k <0\} } \bigr) \Delta t.
\end{equation}
The algorithm relies on the discrete-time Bellman equation from Dynamic Programming Principle (DPP). The Bellman equation between $t_k$ and $t_{k+1}$ is given by:
\begin{equation} \label{eq:Bellman_DP}
             V(t_k,X_k,I_k ) = \inf_{B_k\in \mathcal{B}_k(I_k)}  
\Bigl\{f(X_k,B_k,M_k) \Delta t+ \EE[V(t_{k+ 1} , X_{k+1}, I_{k+1}) \big|\, {X_k,I_k}] \Bigr\}
\end{equation}
where $\mathcal{B}_k(I_k)$ is the discrete-time counterpart of the admissible control set \eqref{constraint_set_t} for $B_k$, adjusted to make sure that $I_{k+1}$ in \eqref{eq:I_discrete} will remain within $[I_{\min}, I_{\max}]$:
\begin{align}\label{eq:Bk}
\mathcal{B}_k(I_k) := \left\{ B_k \in \mathbb{R} : 
B_{\min}(I_k)
\leq B_k \leq 
B_{\max}(I_k)
\right\}
\end{align}
where 
\begin{align*}
B_{\min}(I_k) &:= \max\left( B_{\min}, \frac{\eta \cdot (\Imin - I_k)}{\Delta t} \right), \qquad
B_{\max}(I_k) := \min\left( B_{\max}, \frac{\Imax - I_k}{\eta \Delta t} \right).
\end{align*}
The expectation in \eqref{eq:Bellman_DP} is taken over the random variable $X_{k+1}$, conditioned on current wind power output $X_k$.  In line with the DPP approach, we 
define the cost-to-go $q$-value
\begin{equation}\label{eq:Q}
    Q(t_k,X_k,I_{k+1}):= \EE \Bigl[V(t_{k+1}, X_{k+1}, I_{k+1}) \, \Big| \, X_k\Bigr]
\end{equation}
and characterize the optimal control at step $t_k$ via
\begin{equation} \label{eq:optimal}
 B^*_k(X_k, I_k) := \arginf_{B_k \in \mathcal{B}_k(I_k)} \Bigl\{ f(X_k,B_k, M_k) \Delta t + Q(t_k,X_k,I_{k}+\bigl(\eta B_{k} \mathbb{1}_{ \{B_k  > 0\} } +\frac{1}{\eta} B_{k}\mathbb{1}_{\{B_k <0\} } \bigr) \Delta t) \Bigr\}.
\end{equation}
To solve  \eqref{eq:optimal} we  approximate $Q(t_k, \cdot, \cdot)$ at each $t_k$  using a machine learning model.   

\subsection{SHADOw-GP Algorithm}
The algorithm approximates the $q$-value $Q(t_k, \cdot, \cdot)$ via a statistical emulator $\widehat{Q}_k$ for each $k$ using Gaussian Process Regression \parencite{ludkovski2018simulation}. The emulator is fitted via empirical regression using the training set, termed the simulation design, $\mathcal{D}^v_k$ of size $N$. Note that since the $q$-value  $\widehat{Q}_k$ is a function of the \emph{current} wind power generation $X_k$ at step $t_{k}$ and the \emph{lookahead} SoC $I_{k+1}$ at step $t_{k+1}$, the simulation design is $\mathcal{D}^v_k = (X^n_k, I^n_{k+1})_{n=1}^{N}$.  

Solving \eqref{eq:optimal} requires repeated calls of optimization for the \(N\) training samples in \(\mathcal{D}^v_k\) during training and later for the  post-training simulations. To mitigate this computational cost, we also construct an emulator, \(\widehat{B}_k\), for the control map. Since $\widehat{B}_k$ is a feedback map of \emph{current} wind generation $X_k$ and \emph{current} SoC $I_k$, the respective simulation design is  $\mathcal{D}^b_k = (X^{b,n}_k,I^{b,n}_k)_{n=1}^{N_b}$ of size $N_b$ and distinct from $\mathcal{D}^v_k$. This results in an actor-critic scheme as described in \textcite{hure2020applications}. 

Following the standard DPP procedure, we proceed backward in time, starting with the known terminal condition $\widehat{Q}_{K-1}(X_{K-1},\iota)=g(\iota)$. 
For $k= K-1,\ldots,0$, we repeat the following 6 steps: \medskip

\noindent 1) Generate design spaces $\mathcal{D}^v_{k-1} = (X^n_{k-1}, I^n_{k})_{n=1}^{N}$ and $\mathcal{D}^b_k = (X^{b,n}_k, I^{b,n}_{k})_{n=1}^{N_b}$;\\
2) Evaluate the pointwise optimal control $b_{k}^{*,n}$ for each input $(X_{k}^{b,n},I_{k}^{b,n}), n=1,\ldots, N_b$ in $\mathcal{D}^b_k$ according to \eqref{eq:optimal} using a numerical optimizer;\\
3) Construct the policy map emulator $\widehat{B}_k: (X,I) \mapsto \mathcal{B}_k \subset \mathcal{R}$ by regressing $b_{k}^{*,n}$ against design $\mathcal{D}^b_k$,
i.e., an empirical $L^2$-projection into the given function space $\mathcal{H}^b$,
\begin{equation}
\label{Bfitting}
\widehat{B}_{k}(\cdot)=\arginf_{h_{k} \in \mathcal{H}^b} \sum_{n=1}^{N_b}\left(h_{k}\bigl(X^{b,n}_{k}, I^{b,n}_{k}\bigr)-b^{*,n}_{k}\right)^2;
\end{equation} \\
4) Perform a one-step forward simulation of $\mathcal{D}^v_{k-1}$: $X^n_{k-1} \mapsto  X^n_{k}$ for $n=1,\ldots,N$;\\
5) Evaluate the pointwise value function in \eqref{eq:Bellman_DP} for each $(X^n_k,I^n_k)$ using the control $\hat{b}^n_k = \widehat{B}_k(X^n_k,I^n_k)$:
\begin{equation}
\label{conditional}
        v_{k}^n = f(X^n_{k},\hat{b}^{n}_{k},M_{k} )  \Delta t+   \widehat{Q}_{k}(X^n_{k} , I^n_{k}+\bigl(\eta \hat{b}^{n}_{k} \mathbb{1}_{ \{\hat{b}^{n}_{k}  > 0\} } +\frac{1}{\eta} \hat{b}^{n}_{k}\mathbb{1}_{\{\hat{b}^{n}_{k} <0\} } \bigr) \Delta t), \quad n = 1,\ldots,N;
\end{equation} \\
6) Construct the $q$-value emulator for
$
        \widehat{Q}_{k-1}: (X,I)\mapsto   \mathbb{E} \bigl[\widehat{V}(t_{k}, X_{k}, I) \mid X_{k-1}=X \bigr]
$
by regressing $v_{k}^{1:N}$ against the design $\mathcal{D}^v_{k-1}$ 

\noindent where steps steps 4--6 involving the design $\mathcal{D}^v_{k-1}$ are skipped at final iteration $k=0$. In the subsequent sections, we delve deeper into our selection of simulation design, $q$-value emulator, and control emulator.


\subsection{Simulation Design }
\label{sim_design}
The algorithm requires the bi-variate simulation design $\mathcal{D}^{v}_k$ at each time $t_k$ for $\widehat{Q}_k(\cdot)$. The regression emulator $\widehat{Q}_k(\cdot)$ will be more accurate in regions with many samples and worse in regions with few samples. This necessitates a simulation design with thorough coverage.  The standard space-filling approach is to
take $(X^n_{k},I^n_{k+1}) \sim $Unif$([X_{\min},X_{\max}]\times [I_{\min},I_{\max}])$. However, uniform i.i.d.\ samples tend to cluster, resulting in higher variance of the emulator. 
In this work, we rely on Latin Hypercube Sampling (LHS) to pick our design sites. 
At each timestep $t_k$, the boundaries of the design $\mathcal{D}^k_v$ should account for the support of the random variables  $I^*_{k+1}$ and $X_{k}$. While $X_{k}$ is exogenous, $I^{*}_{k+1}$ is controlled by $B^{*}_k$. 
To this end, we space fill over the training domain $[X^{k}_{\min},X^{k}_{\max}]\times [I_{\min},I_{\max}]$, where the choice of $X_{\min}^{k}$ and $X_{\max}^{k}$ is taken as three times the standard deviation $\STD(X_k)$. We further augment our design with \textbf{replication} to decrease the variance of estimate $\widehat{Q}_k(\cdot)$. i.e; we divide our training design into $N_{loc}$ distinct sites, with each  distinct input repeated $N_{rep}$ times (for the remainder of the subsection,  $\bx \equiv (X,I)$ denotes a generic training input):
\begin{equation}
\mathcal{D}^v_k=\{\underbrace{\bx^{1,1}, \ldots, \bx^{1,N_{rep}}}_{N_{\text {rep }} \text { times }}, \underbrace{\bx^{2,1}, \ldots,\bx^{2,N_{rep}}}_{N_{\text {rep }} \text { times }},\ldots, \underbrace{\bx^{N_{loc},1},\ldots,\bx^{N_{\text {loc }}}}_{N_{\text {rep }} \text { times }}\}.
\end{equation}
\noindent 
The total simulation budget at each step $k$ is thus $|\mathcal{D}^v_k| = N_{loc} \times N_{rep}$. Subsequently, one-step forward simulations and evaluations are performed to acquire the respective responses $v^{1,1}, \ldots, v^{i, j}, \ldots, v^{N_{\mathrm{loc}}, N_{\mathrm{rep}}}$. Denote by 
\begin{equation}
    \bar{\mathcal{D}}^v_k = \{\bx^1,\bx^2,\ldots,\bx^{N_{loc}}\}
\end{equation} 
the design $\mathcal{D}^v_k$ without replicates, i.e., $\bar{\mathcal{D}}^v_k$ contains just the $N_{loc}$ unique training samples. 
After pre-averaging the responses of each replicated batch $ \bar{v}^n:=\frac{1}{N_{\mathrm{rep}}} \sum_{j=1}^{N_{\mathrm{rep}}} v^{n, j}
$,  the regression model for the continuation value emulator $\widehat{Q}_k(\cdot)$ is constructed by regressing $\bar{v}^{1:N_{loc}}$ against the design $\bar{\mathcal{D}}^v_k := (\bx^{1:N_{loc}}$). The pre-averaging is equivalent to a Monte Carlo approximation of conditional expectation in \eqref{eq:Q}. The replicated design lowers regression training errors thanks to the decreased variability in $\bar{v}^n$'s, raising the signal-to-noise ratio. Furthermore, to ensure the existence of samples at the minimum SoC $\Imin$ and the maximum SoC $\Imax$, we add additional  training inputs at $I=I^{k+1}_{\min}$ and $I=I^{k+1}_{\max}$. Similarly, we opt for additional training points at $X=X^k_{\min}$ and $X=X^k_{\max}$. That is, we deploy a \textbf{fencing} mechanism by creating an additional layer of (evenly spaced) training inputs on the boundary of the sampling domain $[X^k_{\min},X^k_{\max}]\times [I^{k+1}_{\min}$  $,I^{k+1}_{\max}]$ prior to replication, see Figure \ref{fig:Training design} below.

We also use LHS to generate the training design \(\mathcal{D}^b_k = (X^{b,n}_k, I^{b,n}_{k})_{n=1}^{N_b}\) for \(\widehat{B}_k(\cdot)\), covering the same domain as \(\mathcal{D}^v_k\). 
Since the optimal control values \eqref{eq:optimal} do not depend on forward simulations, \(\mathcal{D}^b_k\) does not require replicates.

\begin{remark}
Unlike reinforcement learning, our algorithm requires a pre-specified space-filling design over the domain $[X_{\min},X_{\max}] \times [I_{\min},I_{\max}]$, whereas RL relies on exploration to generate such coverage. In practice, our design size is on the order of only a few hundred unique design points, while RL typically requires tens of thousands to millions of samples to achieve comparable accuracy. The drawback, however, is that our approach is hindered by the need for either pre-specified or fitted system dynamics to generate such design apriori, whereas RL can operate in a model-free manner provided sufficiently large number of historical trajectories are available.

\end{remark}

\subsection{Gaussian Process Emulator for \QQ}

In order to enhance the numerical optimization of $B^*_k$ in \eqref{eq:optimal}, we opt for an emulator that has an analytical gradient. To achieve this, we make use of Gaussian Process regression (GPR). GPR models $Q_k(\cdot)$ as a Gaussian Process  specified by a mean function $m(\bx)$ (taken to be zero after standardizing the outputs) and positive definite covariance function  $c(\bx,\bx^{\prime}; \vartheta)$ \parencite{books/lib/RasmussenW06}.  The hyperparameters $\vartheta$ of the covariance kernel $c(\cdot,\cdot)$ specify the smoothness of $\widehat{Q}$.  Given a training design $(\bx^{1:N}, v^{1:N})$ and an input $\bx^*$, the continuation value, $\widehat{Q}(\bx^{*})$, is the posterior mean of the GP given by 
\begin{equation}\label{eq:gp-mean}
\widehat{Q}(\bx_*) 
= \mathcal{C}^\top_*(\bC+ \sigma_{\epsilon}^2 \bI)^{-1} \mathbf{v} 
= \bigl[c(\bx_*, \bx^1;\vartheta), \ldots, c(\bx_*, \bx^N;\vartheta)\bigr]
(\bC+ \sigma_{\epsilon}^2 \bI)^{-1} \mathbf{v},
\end{equation}
where $\bC$ is the $N \times N$ covariance matrix with $\bC_{k,l} = c(\bx_k, \bx_l;\vartheta)$, $\bI$ is $N\times N$ identity matrix,  
$\mathbf{v}=\left[v^1, \ldots, v^N\right]^\top$, and
$\sigma_{\epsilon}^2$ represents observation noise. 
One may opt for different kernels, see \textcite{books/lib/RasmussenW06}. In our study, we choose the twice-differentiable anisotropic Mat\'ern-5/2 kernel with three hyperparameters $\vartheta = (\sigma^2_{p},\ell_{ 1},\ell_{ 2})$
\begin{equation}
\begin{split}
    &c_{M 52}(\bx, \bx^{\prime} ; \vartheta):=\sigma_p^2 \prod_{j=1}^2 \Bigl(1+\frac{\sqrt{5}}{\ell_j}|x_j-x_j^{\prime}|+\frac{5}{3 \ell_j^2}(x_j-x_j^{\prime})^2 \Bigr) \cdot \exp \bigl(-\frac{\sqrt{5}}{\ell_j}|x_j-x_j^{\prime}| \bigr),
\end{split}
\end{equation}
where $\sigma^2_{p}$ indicates the magnitude of the response, and $\ell_{ 1}$ and $\ell_{2}$ determine how the response fluctuates with respect to wind power generation (MW) and SoC (MWh), which are expressed in different scales and units and hence have different lengthscales. In simple terms, a large lengthscale suggests a smooth response surface, while a small lengthscale indicates a non-smooth surface with significant fluctuations. 
The hyperparameters $\vartheta$ and $\sigma_{\epsilon}^2$ are optimized using maximum likelihood estimation (MLE). To accelerate GP $\widehat{Q}_k(\cdot)$ optimization at step $k$, we warm-start with the hyperparameters $\vartheta^{(k+1)}$ obtained from the trained GP $\widehat{Q}_{k+1}(\cdot)$ at step ${k+1}$. We also standardize both inputs and outputs in the range $[-1,1]$ to improve the numerical stability of the GP MLE. 


\subsection{Gaussian Process emulator for \BB}
The constrained optimization problem in \eqref{eq:optimal} is given implicitly in terms of the emulator $\widehat{Q}_k(\cdot)$, with  the respective first-order-condition tied to $\partial \widehat{Q}_k/\partial I$.
GP allows the use of faster gradient-based optimizers thanks to its analytical gradients. Differentiating the GP $\widehat{Q}(\cdot)$ with Mat\'ern-5/2 kernel kernel in $x_j$
gives another GP with  posterior mean at input $\bx_{*}$ given by 
\begin{equation}\label{eq:grad}
\frac{\partial \widehat{Q}}{\partial x_j}(\bx_*) 
= \sum_{n=1}^{N} \alpha_n 
\frac{\partial c_{M52}}{\partial x_j}\bigl(\bx_*, \bx^n ; \vartheta\bigr) 
= \sum_{n=1}^{N} \alpha_n\, c_{M52}\bigl(\bx_*, \bx^n;\vartheta\bigr)
\frac{-\frac{5}{3 \ell_j^2}r_j-\frac{5^{3 / 2}}{3 \ell_j^3}r_j|r_j|}
{1+\frac{\sqrt{5}}{\ell_j}|r_j|+\frac{5}{3 \ell_j^2}r_j^2}.
\end{equation}
where $\alpha_n$ is the $n$-th component of $\left(\bC+\sigma_\epsilon^2 \bI\right)^{-1} \mathbf{v}$ and
with $r_j = (x_j-x_j^{*})$. In our case, we utilize the gradient with respect to SoC, $j=2$ in \eqref{eq:grad}. The presence of the black-box $\widehat{Q}$ makes the optimization problem non-convex. Since we have one-dimensional control, we project the constraints onto the optimal control solution from the unconstrained case.  We employ the unconstrained, gradient-based $\textbf{L-BFGS}$ solver from the \texttt{SciPy} library to retrieve optimized outputs without constraints. That is, instead of solving for \eqref{eq:optimal} with $\hat{Q}_k(\cdot)$ directly, we instead solve:
\begin{equation} \label{eq:optimal_unconstrained}
     B^{\dagger}_k(X_k, I_k) :=\arginf_{B_k \in \mathbb{R}} \Bigl\{ f(X_k,B_k,M_k) \Delta t + \widehat{Q}_k(X_k,I_{k}+\bigl(\eta B_{k} \mathbb{1}_{ \{B_k  > 0\} } +\frac{1}{\eta} B_{k}\mathbb{1}_{\{B_k <0\} } \bigr) \Delta t) \Bigr\}.
\end{equation}
which is then projected back onto $\mathcal{B}_k$. As an alternative, we may find the zero of the first-order condition of \eqref{eq:optimal_unconstrained}.
When the cost function is given by \( f(X_k,B_k,M_k) = (X_k-M_k-B_k)^2 \), the first-order condition for the optimal \( B_k^\dagger \) is:

\begin{equation}\label{eq:foc}
2\big(X_k - M_k - B_k^\dagger\big)
= \frac{d\,\widehat{Q}_k}{dB}(X_k,I_{k}+\bigl(\eta  B_k^\dagger \mathbb{1}_{ \{ B_k^\dagger  > 0\} } +\frac{1}{\eta}  B_k^\dagger\mathbb{1}_{\{ B_k^\dagger <0\} } \bigr) \Delta t)
\end{equation}
Using the analytical gradient of the GP, \( \frac{\partial \widehat{Q}_k}{\partial I} \) and chain-rule, we can apply a root-finding algorithm to \eqref{eq:foc} to determine \( B^\dagger_k \). Then, we directly enforce the constraints that define the feasible set $\mathcal{B}_{k}$ onto $ B^{\dagger}_k(X_k,I_k)$ to obtain the optimal pointwise control ${B}^{*}_k(X_k,I_k)$ in \eqref{eq:optimal}. 

\textbf{Emulating the feedback control map:} The unconstrained optimization with manual projection in \eqref{eq:optimal_unconstrained} serves two purposes: saving training time compared to constrained optimization and construction of the emulator for $\widehat{B}_k(\cdot)$. Directly building an emulator over the output ${B}^{*}_k(X_k,I_k)$ with hard constraints does not guarantee feasibility. 
To guarantee feasibility, we train $\widehat{B}_k$ by regressing \textit{unconstrained} optimized $(b^{\dagger, n}_{k})_{n=1}^{N_b}$ against $\mathcal{D}^b_k = (X^{b,n}_{k}, I^{b,n}_{k})_{n=1}^{N_b}$ with a projection layer onto feasible set $\mathcal{B}_k$ post-training. The equation for $\widehat{B}_k$ \eqref{Bfitting} becomes
\begin{equation}
\label{eq:Bhat_noproject}
\widecheck{B}_k(\cdot)= \arginf_{h_{k} \in \mathcal{H}^b} \sum_{n=1}^{N_b}\left(h_{k}\left(X^{b,n}_{k}, I^{b,n}_{k}\right)-b^{\dagger,n}_{k}\right)^2;
\end{equation} 
\begin{equation}
\label{eq:Bhat_project}
    \widehat{B}_k(X^{b,n}_k, I^{b,n}_k) = \operatorname{Proj}_{\mathcal{B}_k(I^{b,n}_k)}  \{\widecheck{B}_k(X^{b,n}_k, I^{b,n}_k)\},
\end{equation}
where the projection depends pointwise on $I^{b,n}_k$. To fit \( \widecheck{B}_k \), we rely on GP with a Matérn-\(3/2\) kernel. We experimented with three different kernels—Exponential, Matérn-\(5/2\), and Matérn-\(3/2\), and found that the Matérn-\(3/2\) kernel provided the best performance. Our approach guarantees feasibility of $\widehat{B}_k(X,I)$;  see \textcite{ RL_1dprojection} for similar projection of BESS within RL framework. The complete SHADOw-GP procedure is outlined in Algorithm~\ref{alg:RMCalgorithm}. Recall that $\bar{\mathcal{D}}^v_k$ is the design $\mathcal{D}^v_k$ without replicates. Here, \(\mathcal{N}_{loc} = \{1,2,\dots,N_{loc}\}\), \(\mathcal{N}_{rep} = \{1,2,\dots,N_{rep}\}\), and 
\(\mathcal{N}_{b} = \{1,2,\dots,N_{b}\}\).
 
\textbf{Estimating the value function: }
After training, Algorithm~\ref{alg:RMCalgorithm} yields continuation-value emulators $\{\widehat{Q}_k(\cdot)\}_{k=1}^{K-1}$ and optimal control emulators $\{\widehat{B}_k(\cdot)\}_{k=0}^{K-1}$. To evaluate the resulting hybrid resource output trajectory $(O_k)$ and the respective value function, we utilize Monte Carlo simulation.  Given an initial state $(X_0,I_0)$, we generate $M$ out-of-sample paths $(X^m_{0:K}, {\hat{I}}^{m}_{0:K})$, $m = 1,\ldots,M$, where the optimized SoC ${\hat{I}}^{m}_{k+1}$ is based on $\hat{b}^{m}_k:=\widehat{B}_k( X^m_k,{\hat{I}}^{m}_k)$  evaluated from \eqref{eq:Bhat_project}. This gives cumulative pathwise realized cost 
\begin{equation}\label{eq:pathwise}
    v^m_{0:K} = \sum_{k=0}^{K-1} f(X^m_k,\hat{b}^{m}_k,M_k) + g({\hat{I}}^{m}_T),  
\end{equation} and the resulting Monte Carlo estimate of the value function: 
\begin{equation}\label{eq:value}
    \widehat{V}({0,X_0,I_0}) = \frac{1}{M}\sum_{m=1}^{M} v^m_{0:K}.
\end{equation}

\begin{algorithm}[ht]
    \caption{Stochastic hybrid asset dispatch optimization with Gaussian Process (SHADOw-GP) }
    \label{alg:RMCalgorithm}
    \begin{algorithmic}[1]
        \STATE  $\textbf{Input:}$ $K$ steps, $N_{loc}$ sites, $N_{rep}$ replications, $N_{b}$ numerical optimizations per step.
        \STATE Set $\widehat{Q}_{K-1}(X,I)=g(I) $ (No emulation)

        \FOR{k = $K-1$ to $0$}
            \STATE Generate training design $\mathcal{D}^b_{k} = \{(X^{b,i}_{k},I^{b,i}_{k}),  i\in \mathcal{N}_{b}\}$.
            \STATE Optimize $b^{\dagger,i}_{k}$ in (\ref{eq:optimal_unconstrained}) pointwise for each $(X^{b,i}_{k}, I^{b,i}_{k})$ for $i\in \mathcal{N}_{b}. $ 
    
            \STATE Fit control GP $\widecheck{B}_{k}(\cdot) $  according to \eqref{eq:Bhat_noproject} using $b^{\dagger,i}_{k}$ and $(X^{b,i}_{k}, I^{b,i}_{k})$  for $i \in \mathcal{N}_{b}. $ 
        \IF{$k \neq 0$}
            \STATE Generate replicated training design $\mathcal{D}^v_{k-1} = \{(X^{i,j}_{k-1},I^{i,j}_{k}), i\in \mathcal{N}_{loc} \text{ and }j \in \mathcal{N}_{rep}\}$.
            \STATE Generate one-step paths: $X^{i,j}_{k-1} \mapsto X^{i,j}_{k}$ for $i\in \mathcal{N}_{loc}$ and $j \in \mathcal{N}_{rep}$.
            \STATE Evaluate $\hat{b}^{i,j}_{k}$ according to \eqref{eq:Bhat_project} using $\widecheck{B}_k(\cdot)$ on $(X^{i,j}_{k},I^{i,j}_{k})$ for $i\in \mathcal{N}_{loc}$ and $j \in \mathcal{N}_{rep}$.
            \STATE Evaluate $v^{i,j}_{k}$ in (\ref{conditional}) for $i\in \mathcal{N}_{loc}$ and $j \in \mathcal{N}_{rep}$.
            \STATE Average over replicates: $\bar{v}_{k}^{n}=\frac{1}{N_{rep}} \sum_{j=1}^{N_{rep}} v_{k}^{n, j}$ for $n\in \mathcal{N}_{loc}$
            \STATE Fit continuation-value GP $\widehat{Q}_{k-1}(\cdot)$ by regressing $\bar{v}_{k}^{n}$ against $(X^{n}_{k-1}, I^{n}_{k})$ in reduced design $\bar{\mathcal{D}}^v_{k-1}$ for $n\in \mathcal{N}_{loc}.$ 
            \ENDIF

        \ENDFOR
        
    \end{algorithmic}
    $\textbf{Output:}$  $q$-value emulators $\{\widehat{Q}_k(\cdot)\}_{k=1}^{K-1}$ and control emulators $\{\widehat{B}_k(\cdot)\}_{k=0}^{K-1}$
\end{algorithm}

\subsection{Toy Linear-Quadratic (LQ) Case Study with Stationary Dynamics}\label{ssec:4.5}
We begin with a toy LQ example under stationary dynamics, for which an analytical solution is available. This experiment benchmarks the algorithm against the LQ solution. Specifically, we use the analytical solution from Proposition \ref{lq_control} as a baseline for evaluating the SHADOw-GP algorithm. We consider quarter-hour increments $\Delta t = 1/4$ with time-stationary generation with a constant dispatch target, $M_k \equiv 5 $ MW for $k=0,\ldots,95$. 
The parameters of the stationary Jacobi process and BESS are listed in Table \ref{tab: stationary table}. Minimum generation is set to $X_{\text{min}} = 0$ and maximum generation is $X_{\text{max}} = 10$MW. Our choice of $B_{\max}$ and $B_{\min}$ implies that the charge / discharge capacity can handle deviations up to $1\STD(X_k) = 1$MW away from the dispatch target of $5$MW. We take a  $3$-hour battery with $I_{\max}= 3 B_{\max} = 3$MWh \parencite{CAISO_BESS}. In our simulations, the BESS starts with $50\%$ SoC, $I_0 = 1.5$MW; accordingly the terminal condition is $g(\iota)=\mathcal{P}(\iota-I_{\text{target}})^2$ with $\mathcal{P}=10$ and $I_{\text{target}} = 1.5$MW.

\begin{table}[ht]
\centering
\begin{tabular}{c|c|c}
\hline
$\alpha_k= 0.5$  & $m_k= 5$ (MW)&  $\sigma_k = 0.2 \ (hr^{-1})$   \\
\hline
$\text{SoC}_{\min}= 0$ & $\text{SoC}_{\max}= 1$ & $\text{I}_{\text{max}}= 3 
 \ (\text{MWh})$
\\
\hline
$B_{\text{min}} = -1 (\text{MW})$ & $B_{\text{max}} = 1 \ (\text{MW})$ &  $\eta = 1$  \\
\hline
\end{tabular}
\caption{Parameters for the stationary example in Sect. \ref{ssec:4.5}}\label{tab: stationary table}
\end{table}

The GP emulators $\{\widehat{B}_k(\cdot),\widehat{Q}_k(\cdot)\}_{k=0}^{95}$ are trained via LHS designs on $\mathcal{D}^v_k$ with size $ N = N_{loc}\times N_{rep} = 640 \times 50$ and $\mathcal{D}^b_k$ with $N_b = N_{loc}$ respectively, at each time step $k$. Here, $40$ of the $N_{loc}$ design points come from the fencing of the boundary of sampling domain. The implementation of SHADOw-GP uses the Python \texttt{scikit} library run on a MacBook Pro M1 laptop with 16GB RAM and an Apple M1 8-core CPU. The computation takes approximately 15 minutes.

\begin{figure}[!ht]
\begin{center}
    \includegraphics[width=0.40\linewidth]{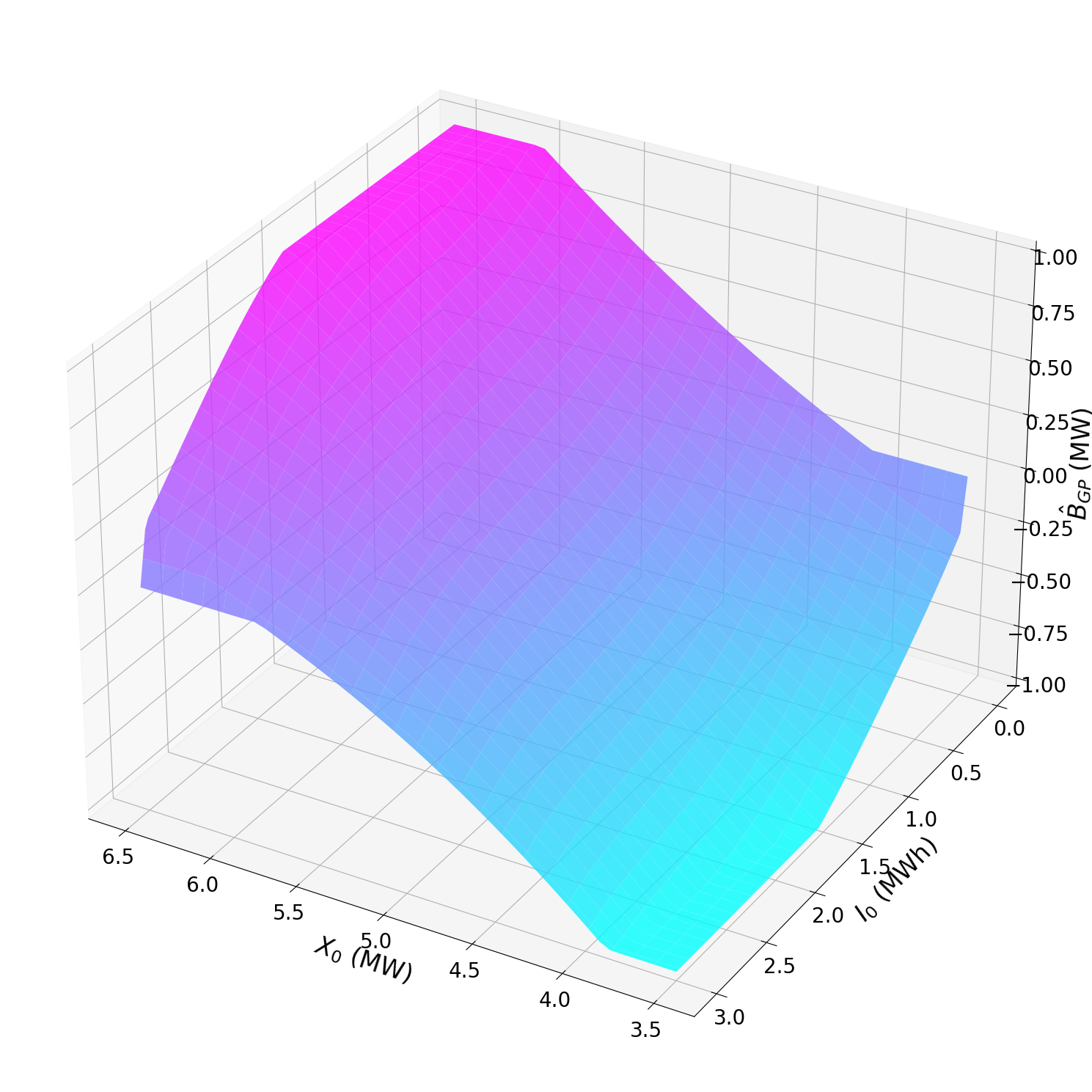}
    $\qquad$ \includegraphics[width=0.47\linewidth]{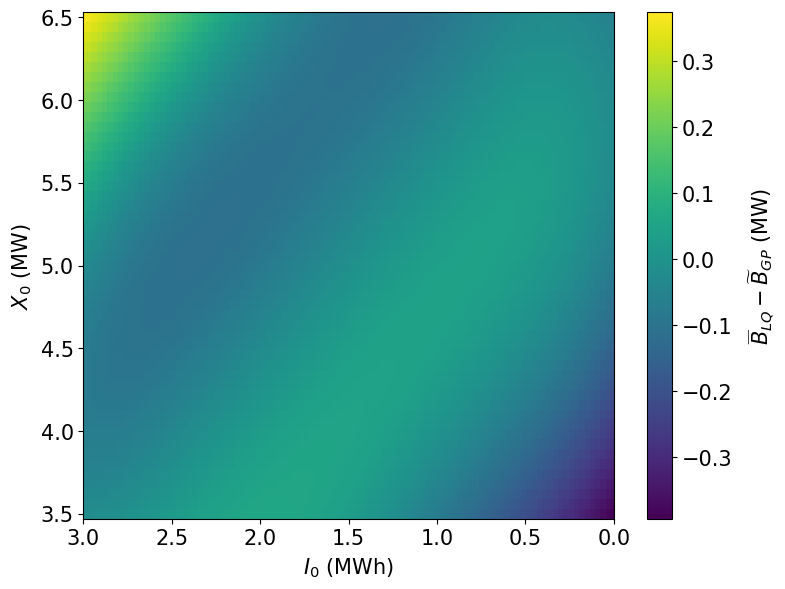}
    \caption{ \emph{Left panel:} Learned control policy $(X,I) \mapsto \widehat{B}_{GP,k}(X,I)$ for the setting of Section \ref{ssec:4.5} at $k=0$ as a function of SoC $I$ and wind generation $X$. The dispatch target is $M_0=5$MW. \emph{Right panel:} Difference  between the LQ and GP-based controls without power and capacity constraints, $\bar{B}_{LQ,0}(X,I)- \widecheck{B}_{GP,0}(X,I)$ at $k=0$. }
    \label{fig:GPvLQsurface}
\end{center}
\end{figure}

The left panel of Figure \ref{fig:GPvLQsurface} visualizes the resulting policy $\widehat{B}_{GP,k}(\cdot, \cdot)$ at $k=0$, based on the quadratic running cost \eqref{eq:tildeF}. Observe that when SoC is far from empty or full, $\widehat{B}_{GP,0}(X,I) \simeq X-M_0$ is almost linear in the middle of the policy surface. However, when both the SoC $I$ and renewable generation $X$ are high, the optimal control decreases the charging rate $\widehat{B}_{GP,0}(X,I) < X-M_0$ to maintain SoC headroom. Similarly, the controller throttles the BESS discharging when the SoC is very low. As a result, over the day $\hat{I}_k$ tends to stay in a ``safe zone'' and away from SoC limits, demonstrating precautionary risk mitigation behavior. 

Next, we compare our SHADOw-GP algorithm with the closed-form LQ solution of Prop.~\ref{lq_control} which approximates control and state constraints. This yields strategy $\bar{B}_{LQ,k}(X,I)$.  
Thanks to the symmetry of the terminal cost with our quadratic penalization on SoC in \eqref{unconstrained_HJB} and $M_k=m_k$, our optimal control in Prop.~\ref{lq_control} simplifies to

\begin{equation}\label{eq:B(c1,c2)}
    \bar{B}_{LQ,k}(X,I; c_1,c_2) = -\kappa P_1(t_k) \cdot ({I} - I_m) + \frac{\kappa}{2} (2-P_2(t_k)) \cdot (X-m). 
\end{equation}
We compare the SHADOw-GP algorithm against the optimal pair \((c_1, c_2)\) for the LQ control policy \eqref{eq:B(c1,c2)}. Specifically, we evaluate MC estimate of the value function \eqref{eq:value} using \(\bar{B}_{\mathrm{LQ},k}(\cdot; c_1, c_2)\), manually projected onto \(\mathcal{B}_k\), over 10,000 Monte Carlo trajectories. This procedure is repeated on a parameter grid \([0,10] \times [0,10]\) with increments \(\Delta c_1 = \Delta c_2 = 0.01\). We find that the pair \(c_1^{*} = 0.08\) and \(c_2^{*} = 0.06\) achieves the lowest cumulative cost, which is approximately 4--6\% higher than that of SHADOw-GP. The difference in optimal controls between SHADOw-GP and LQ with these $c^{*}_1,c^{*}_2$ is shown in the right panel of Figure \ref{fig:GPvLQsurface}. 
The SHADOw-GP charges/discharges less than the LQ strategy with increasing/decreasing SoC and wind power output. In addition to the quadratic ``soft'' approximation of the constraints, this discrepancy also stems from applying the continuous-time control \(\bar{B}_{LQ,k}(\cdot)\) within a discrete-time framework.

\begin{figure}[!ht]
\begin{center}
      \includegraphics[width=0.6\linewidth]{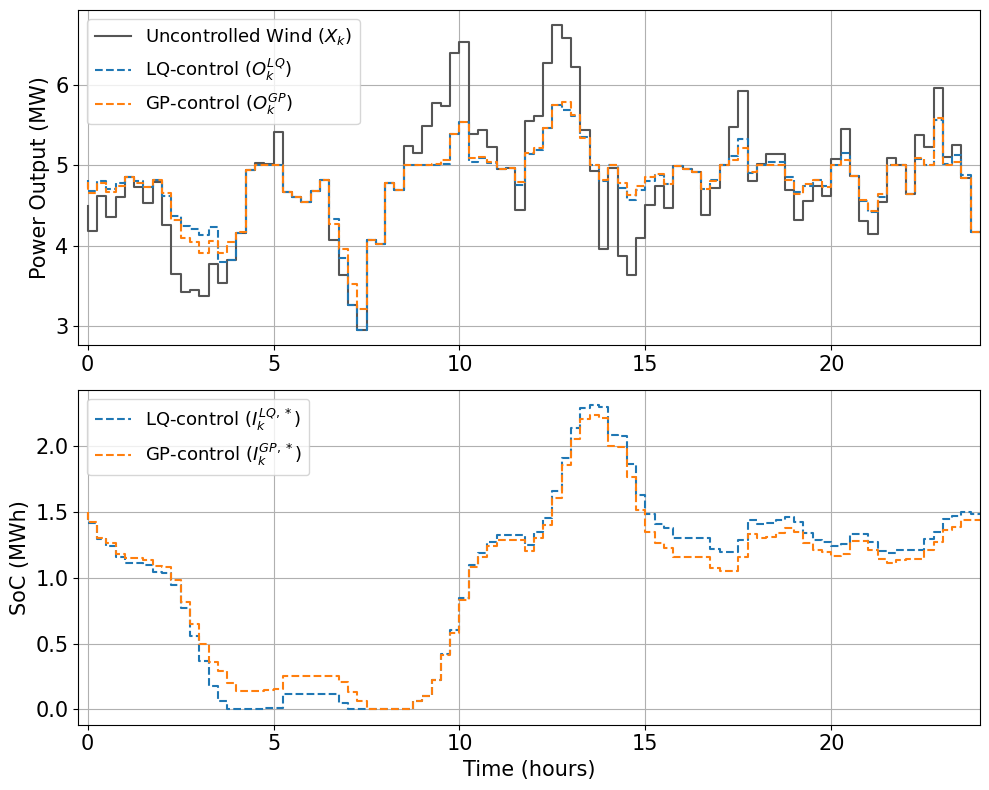}
  \caption{\emph{Top panel:} A trajectory of $(X_k)$ following \eqref{eq:OU} with constant mean $\EE[X_k]=5$ along with firmed hybrid output $(O_k)$ following LQ and SHADOw-GP controls. \emph{Bottom panel}: Corresponding SoC trajectories $({I}^{LQ,*}_k)$ and $({I}^{GP,*}_k)$. }
  \label{fig:Atrajectory}
\end{center}
\end{figure}

The top panel of Figure \ref{fig:Atrajectory} displays a representative 24-hour trajectory of the raw renewable generation $(X_k)$, relative to its firmed  $(O^{LQ}_k)$ and  $(O^{GP}_k)$ outputs based on the analytical LQ and SHADOw-GP control maps.
Both controllers act strategically, conserving some battery capacity rather than myopically charging/discharging to the maximum extent. As expected, larger deviations are adjusted more  aggressively. The GP-based control is better at avoiding SoC limits, for instance keeping ${I}^{GP,*}_k > 0$ non-empty during hours 3-7 (bottom panel of Figure \ref{fig:Atrajectory}) which in turn leads to a smoother overall  trajectory of $(O^{GP}_k)$ compared to $(O^{LQ}_k)$. 

\section{Experiments on Texas-7k grid}\label{sec:5}

In this section, we move beyond the toy LQ setting and evaluate the performance of our SHADOw-GP algorithm under realistic grid conditions. Our experimental platform utilizes a synthetic model of the ERCOT grid, known as the Texas-7k \parencite{TexasGrid}, comprised of 185 renewable assets comprised of 36 solar and 149 wind units, Figure~S1 in the supplementary materials. None of the Texas-7k assets are hybrid as the testbed was developed before the widespread adoption of BESS. Thus, we frame our study as impact analysis for a potential retrofit of existing wind farms with BESS; such retrofits are popular since they face simplified interconnection queue permitting and moreover enjoy tax breaks, enhancing their business case over new greenfield hybrid projects. 

To apply SHADOw-GP to Texas-7k, we first calibrate the dynamics of $(X_k)$ to each individual unit, as described in Section \ref{ssec:5.1}. In Section \ref{ssec:5.2}, we evaluate SHADOw-GP in a synthetic model-based setting by retrofitting the wind farms into hybrid assets. The experiments are conducted with a 3-hour battery and a maximum charge/discharge rate set to 10\% of the unit’s nameplate capacity.
Finally, we assess the benefits of hybridization by feeding the optimized hybrid dispatch \((O_k^*)\) into a full unit commitment and economic dispatch (ED) framework, as described in Section \ref{ssec:4.3}. Specifically, we compare ED costs between the original wind asset and the wind asset retrofitted with a battery, thereby quantifying the system-level impact of hybridization. To this end, we utilize \texttt{Vatic} ED optimization software suite \parencite{Vatic}. Here, we use respective day-ahead forecasts as our dispatch target $(M_k)$ and compute $(O_k^*)$ along actual trajectories of wind generation, rather than on model-based test simulations.

\subsection{Calibration of Wind Dynamics}\label{ssec:5.1}
We utilize the reanalyzed NREL data set \parencite{xu2017creation} for calibration of wind assets. The data set comprises hourly ($\Delta t = 1$) actual and forecasted (as of noon on the day prior) wind generation profiles in 2018 for the 149 wind units in Texas-7k, with generation capacities ranging from 30MW to 352MW. 
Rescaling the raw generation of each wind farm by its nameplate capacity yields \textit{generation ratios} taking values in $[0,1]$. To distinguish from our model-based quantities, we denote by $A^{\ell,d}_k$ the observed generation ratio of asset $\ell$ on day $d$ in hour $k$ and by $F^{\ell,d}_k$ the corresponding day-ahead forecast, later also used as dispatch target $M^{\ell,d}_k$. 

Wind generation is driven by wind speed which exhibits diurnal variations, marked by fluctuations in both mean and variance during a 24-hour period. 
To incorporate these features, rather than using the discrete-time counter part of Jacobi SDE \eqref{eq:OU}, we found that a more flexible fit is provided by 
taking both the drift and diffusion coefficients of $(X_k)$ to be forecast- and time-of-day dependent, with the following  mean-reverting dynamics: 
\begin{equation}\label{eq:discrete calibrated OU}
    X^{\ell,d}_{k+1} = X^{\ell,d}_{k} + \alpha(F^{\ell,d}_k) (F^{\ell,d}_k- X^{\ell,d}_{k}) \Delta t + \sigma(F^{\ell,d}_k) \cdot \epsilon_k \qquad\text{ for } k = 0,1,\ldots,23,
\end{equation}
\noindent 
where $(\epsilon_k)$ is the exogenous noise. Direct calibration of
 \(\alpha(\cdot)\) and \(\sigma(\cdot)\) as continuous functions of the day-ahead forecast ratio \(F \in [0,1]\) in \eqref{eq:discrete calibrated OU} leads to a complex statistical optimization problem  that falls beyond the scope of this paper. To simplify, we take  \(\alpha(\cdot)\) and \(\sigma(\cdot)\) to be piecewise constant in $F$.  Specifically, we partition the forecast values into 10 equi-probable bins and estimate \(\alpha\) and \(\sigma\) independently within each bin. We first assign hourly forecasts to the bins using the non-decreasing map
  $R: F \subset [0,1]\mapsto \mathbb{B}:=\{1,2,\ldots,10\}$, denoting by $\mathcal{N}^{\ell}_r = \{(d,k) : R(F^{\ell,d}_k) = r\}$ the day-hour pairs that fall into bin $r$ for asset $\ell$. Since we have hourly observations for $365$ days, we have $N_{\ell}$ = $8,760$ total observations and  $N^{\ell}_r = \lfloor N_{\ell} / 10 \rfloor = 876 $ for all $\ell$  observations in each bin $r \in \mathbb{B}$. 
  
  We then estimate the mean-reversion rate $\alpha^{\ell}_{r}$ and volatility $\sigma^{\ell}_{r}$ in each bin $r\in \mathbb{B}$ and each asset $\ell$ via least squares under the dynamics of \eqref{eq:discrete calibrated OU}, following the  procedure in  \textcite{Heyman_calibrate,BALATA2021640}. The bin mean-reversion rate is based on a linear regression of the next-hour increments of actual generation $(A^{\ell,d}_{k+1} - A^{\ell,d}_k)$ against the current deviation 
  $(F^{\ell,d}_k - A^{\ell,d}_k)$,
$$ 
\alpha^\ell_{r} = \arg \min_{\alpha \in \mathbb{R}} \sum_{(d,k) \in \mathcal{N}^{\ell}_r} \left( (A^{\ell,d}_{k+1} - A^{\ell,d}_k) - \alpha \cdot (F^{\ell,d}_k - A^{\ell,d}_k) \right)^2;
$$ 
and the bin volatility is the empirical standard deviation of the resulting residuals
$$
  \sigma_{r}^\ell = \STD(\mathcal{E}_{r}^\ell)
                \quad \text{where} \quad  
                \mathcal{E}_{r}^\ell := \{ (A^{\ell,d}_{k+1} - A^{\ell,d}_k) - \alpha^{\ell}_{r} (F^{\ell,d}_k - A^{\ell,d}_k) : (d,k) \in N_{r}^\ell \}.$$

\begin{figure}[!ht]
\centering
  \includegraphics[width=0.475\textwidth,trim=0.4in 0.2in 0.0in 0.15in]{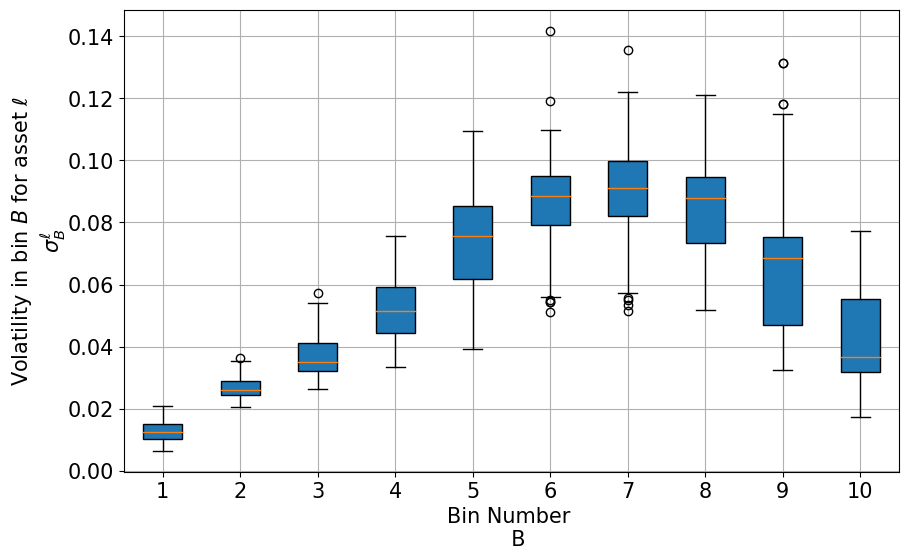}
  \includegraphics[width=0.475\textwidth, trim=0.0in 0.2in 0.4in 0.2in]{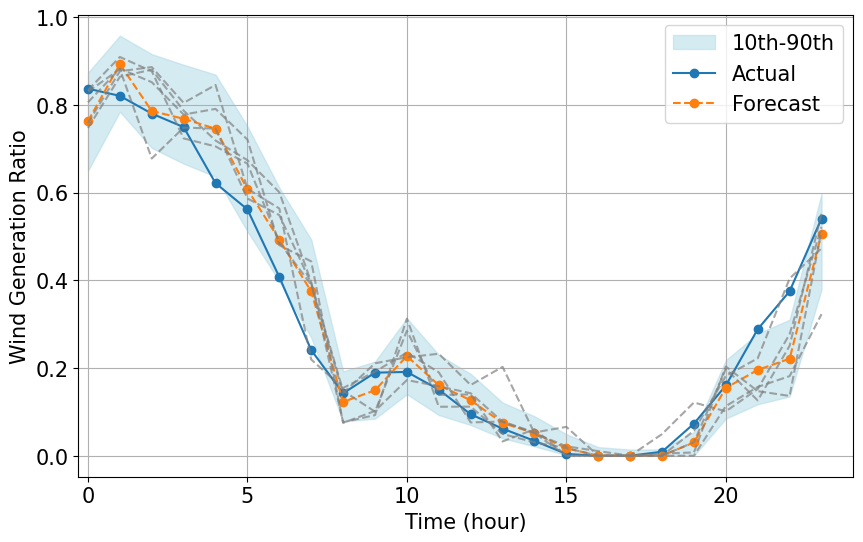}

\caption{
\emph{Left panel}: Boxplot of the calibrated volatility $\sigma^{\ell}_{r}$ across 149 wind assets in Texas-7k, shown as a function of bin $r=1,\ldots, 10$. \emph{Right panel}: Wind generation ratio scenarios generated by \eqref{eq:non-gaussian OU} for Foard City Wind Farm on 2018-04-05. The dotted lines are 5 sample hourly trajectories; the light blue band is the 80\% scenario band. We also show the day-ahead forecast (orange) and the actual generation ratio on that day (blue).} 
\label{fig:Binned Calibration}
\end{figure}

Figure \ref{fig:Binned Calibration} depicts a boxplot of the resulting $\sigma^{\ell}_{r}$ for $r=1,\ldots, 10$ across $\ell=1,\ldots, 149$ assets. The general pattern of $\sigma^{\ell}_{r}$ is low volatility when the forecasted generation ratio is very low or very high (a calm day with minimal wind or strong steady wind), and much higher $\sigma_{r}^\ell$ for generation ratios in the middle. This can be further linked to the nonlinear (approximately cubic) behavior of wind generation as a function of wind speed at intermediate speeds. Figure~S2 in the supplementary materials shows the corresponding boxplot of $\alpha_{r}^\ell$'s. 

It remains to specify the exogenous i.i.d.~shocks $\epsilon_k$ in the calibrated dynamics:
\begin{equation}\label{eq:non-gaussian OU}
    X^{\ell,d}_{k+1} = X^{\ell,d}_{k} + \alpha^{\ell}_{R(F^{\ell,d}_k)} (F^{\ell,d}_k- X^{\ell,d}_{k})  + \sigma^{\ell}_{R(F^{\ell,d}_k)} \cdot \epsilon^{\ell}_{k} \qquad\text{ for } k = 0,1,\ldots,23,
\end{equation}

After calibrating $\alpha^\ell_{\cdot}$ and $\sigma^\ell_{\cdot}$, the binned residuals $\mathcal{E}^\ell_\cdot$ reveal that the shocks $\epsilon^\ell_k$ are non-Gaussian (Supplementary Figure~S3). Moreover, we must enforce $X^{\ell,d}_{k+1} \in [0,1]$ while allowing positive probability of boundary hits, which occur frequently when forecasts satisfy $F^{\ell,d}_k \in \{0,1\}$; see the actual generation ratio trajectories $(A^{\ell,d}_k)$ in Figure~\ref{fig:Binned Calibration}. These features motivate two adjustments to the shock distribution.

\textbf{First}, we replace the Gaussian assumption with a binned bootstrap:
\[
\epsilon^\ell_k \equiv \epsilon^\ell_{R(F^{\ell,d}_k)},
\]
where $\epsilon^\ell_{r}$ is resampled from the empirical residual distribution $\mathcal{E}^\ell_{r}$ for bin $r$.

\textbf{Second}, to ensure boundary hits at $0$ occur with the correct frequency, we modify $\epsilon^\ell_{1}$ (the lowest-decile bin, where $F^{\ell,d}_k=0$) using a mixed distribution with a point mass at zero. We draw $Z^\ell_1 \sim 1-\Bern(p^\ell_1)$ and set
\[
\epsilon^\ell_{1} = \epsilon^{\ell,+}_{1}\,\mathbb{1}_{\{ Z^\ell_1 = 1\}},
\]
where $\epsilon^{\ell,+}_{1}$ is a positive sample from $\mathcal{E}^\ell_{1}$. We estimate $p^\ell_1$ by
\begin{equation}\label{eq:bootstrap}
    \hat{p}^{\ell}_{1} 
    = \frac{\sum_{(d,k) \in \mathcal{N}^{\ell}_{1}} 
      \mathbb{1}_{\{F^{\ell,d}_k=0, F^{\ell,d}_{k+1}=0\}}}{
      \sum_{(d,k) \in \mathcal{N}^{\ell}_{1}} 
      \mathbb{1}_{\{F^{\ell,d}_k=0\}}}.
\end{equation}
Across 149 Texas-7k wind units, $\hat{p}^{\ell}_{1}$ ranges from $72\%$ to $78\%$, ensuring that simulated trajectories hit zero with high probability when $F_k=0$.
Therefore, with high probability, our trajectories hit the boundaries when $F_k = 0$. We apply similar modification of our shocks $ \epsilon^{\ell}_{10}$ where $F^{\ell,d}_k=1$. In this case, the probability of maximum generation $\hat{p}^\ell_{10}$ ranges from $24\%$ to $47\%$. See \textcite{li2024probabilisticspatiotemporalmodelingdayahead} for related multi-site calibration. The right panel of Figure \ref{fig:Binned Calibration} illustrates the behavior of the resulting generation ratios $(X^{\ell,d}_k)$ for a representative wind asset. The shown scenario band (in blue) is based on  $10,000$ Monte Carlo trajectories  $\mathbf{X}^{\ell,d}_\cdot$ generated according to \eqref{eq:non-gaussian OU}, relative to the forecast generation ratio $(F^{\ell,d}_k)$ (in orange). Note the strong tightening of the band (i.e.~minimal variance of $X^{\ell,d}_k$) between hours 15-19 when the forecast is zero, illustrating \eqref{eq:bootstrap}. Supplementary Figure~S4 presents statistical analysis across all the Texas-7k assets to validate \eqref{eq:non-gaussian OU} using an empirical coverage statistic. These results suggest that that the calibrated discretized SDE model performs reasonably well but exhibits some site-to-site variability.

\subsection{SHADOw-GP Evaluation}\label{ssec:5.2}

To assess the impact of hybridization, we consider retrofitting the wind assets in Texas-7k with a three-hour BESS, $\Imax = 3 B_{\max}$ with capacity \( B_{\text{max}} = 0.10 \) representing \( 10\% \) of the nameplate generation capacity. 
Our \( B_{\max} \) is a bit above the variability of the hourly forecast errors, $\mathrm{StDev}( A^{\ell,d}_k - F^{\ell,d}_k) \in [0.05, 0.08]$. 
We take charging efficiency \( \eta = 0.95 \) (10\% round-trip losses; \textcite{LBNL} report that typical roundtrip efficiencies in US are  $75-95\%$) and \( \text{SoC}_{\min} = 0.05, \text{SoC}_{\max} = 0.95 \) (\( 5\% \) SoC buffer).  We use the same simulation design as in the case study in Section \ref{ssec:4.5} with hourly $K=24$ steps. We utilize the same running cost $\tilde{f}$ and terminal cost $g$ as in the stationary example, with terminal penalty $\mathcal{P}=1$. 

\begin{figure}[!htb]
    \centering
        \includegraphics[width=0.45\textwidth, trim={0.2cm 0.95cm 0.cm 0.2cm}]{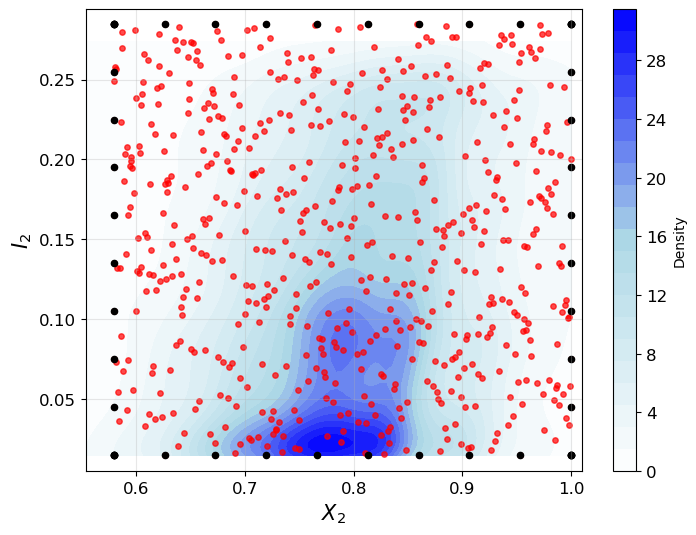}
        \includegraphics[width=0.45\textwidth, trim={0.2cm 0.95cm 0cm 0.2cm}]{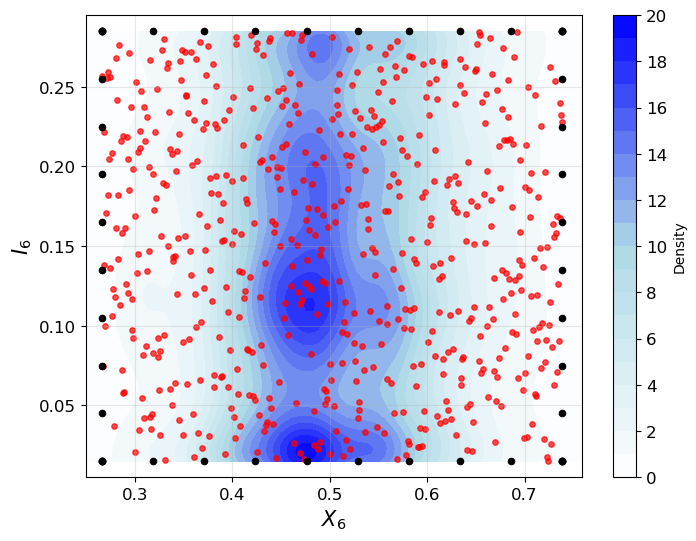}
    \caption{Simulation design $\bar{\mathcal{D}}^v_k$ of size $N_{loc}=600$ for the model calibrated to Foard City Wind Farm. We use LHS on the indicated time-dependent adaptive rectangular training domain; the 40 black dots represent the fencing mechanism. The colors indicate the (kernel-based) density of resulting optimized trajectories  $(X_k,I^*_k)$.
    \emph{Left panel: }  $k=2$. \emph{Right panel: } $k=6$.}\label{fig:Training design}
\end{figure}

The pervasive time-dependence of $(X_k)$ makes the solution, i.e.~the maps $\hat{Q}_k, \hat{B}_k$, also strongly dependent on step $k$.
The SHADOw-GP algorithm automatically adapts the training domain $(\mathcal{D}^v_k)$ at each step $k$ to reflect the time-varying distribution of $(X_k)$, see the $x$-axes in Figure \ref{fig:Training design} that shows $\bar{\mathcal{D}}^v_k$ at two different hours. Figure \ref{fig:Training design} also illustrates the (kernel-smoothed) density of $10,000$ Monte Carlo trajectories of optimally dispatched $(X_k, I^*_k) $ at $k=2$ and $k=6$. The supplementary Figure~S5 plots the differences between the respective control maps at those two hours.

\begin{figure}[!ht]
    \centering
    \includegraphics[width=0.65\linewidth]{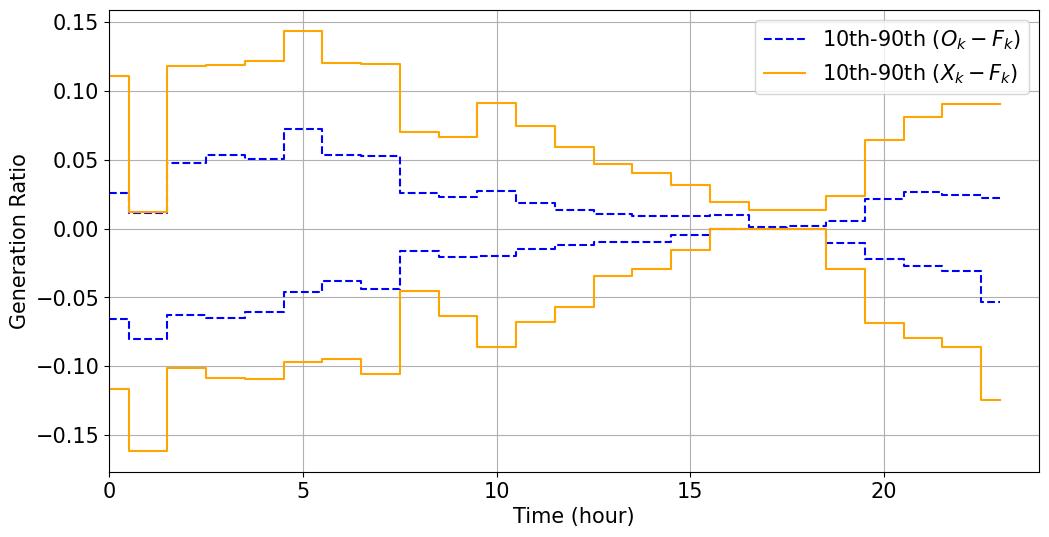}
    \caption{Standalone and hybridized dispatch deviations for Foard City Wind Farm. The $y$-axis is centered around the dispatch target, so perfect match corresponds to zero deviation. Orange solid/blue dotted lines represents the $80\%$ confidence interval of dispatch deviation without BESS $(X_k-F_k)$ and after optimal BESS firming $(O_k-F_k)$. The confidence intervals are based on $10,000$ Monte Carlo trajectories with SHADOw-GP control.}
    \label{fig:bands}
\end{figure}

Finally, Figure \ref{fig:bands} compares the 80\% confidence interval of dispatch deviation without BESS $(X_k-F_k)$ and with BESS $(O_k-F_k)$. As expected, with BESS, we have a tighter band around the target. Notice that the distance between the outer band and the inner band is always less than $B_{\max}$ due to the SoC constraints and the objective $\tilde{f}$. 

\begin{figure}[!ht]
\begin{center}
  \includegraphics[width=0.45\textwidth,trim=0.2in 0.2in 0.2in 0.2in]
  {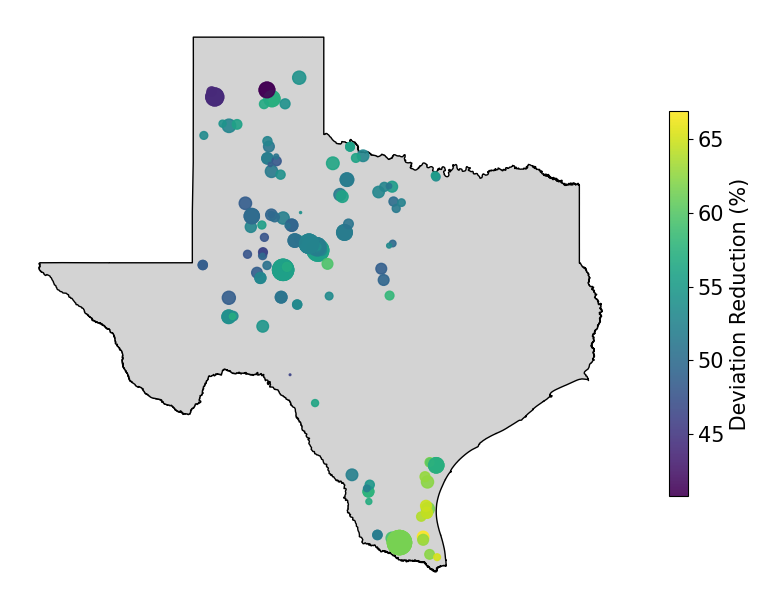}
     \includegraphics[width=0.45\textwidth,trim=0in 0.2in 0.2in 0.2in]{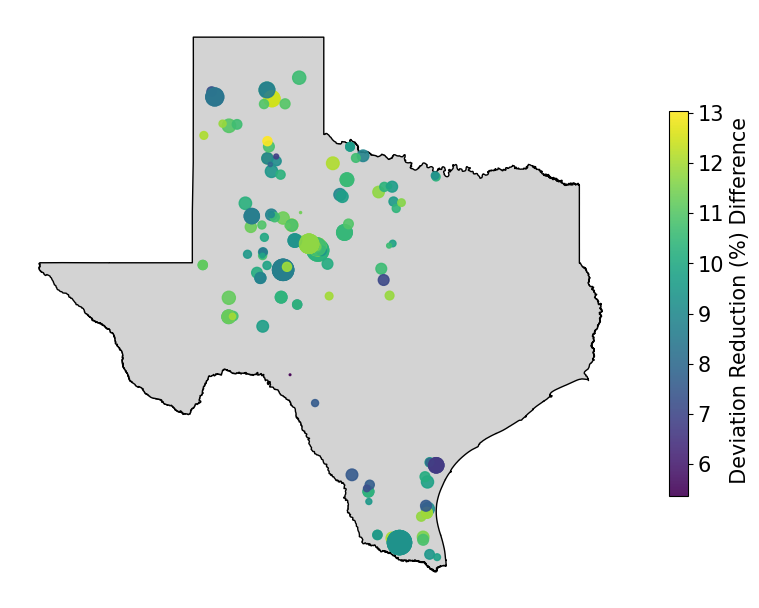}
  \caption{ \emph{Left panel}: Average percentage deviation reduction $\overline{\mathrm{DR}}^{\ell}$ from hybridization. \emph{Right panel}: Improved firming performance (measured as additional reduction in $\overline{\mathrm{DR}}^{\ell}$) for an individual wind asset hybrid retrofit using a  6-hour BESS vis-a-vis a 3-hour BESS. We show the 149 wind assets in the Texas-7k grid, color-coded by average percentage reduction in DR across 24 test days. Symbol size is proportional to nameplate generation capacity. 
 \label{fig:L1reduction}}
\end{center}
\end{figure}

Next, we analyze the results of SHADOw-GP on all wind assets. To do so, we introduce the Deviation Reduction ($\mathrm{DR}$) metric in units (\%) for a given asset $\ell$:
\begin{equation}
\label{eq:DR}
\mathrm{DR}^{\ell,d}
 := \left( \frac{\Dev{A^{\ell,d}} - \Dev{O^{*,\ell,d}}}{\Dev{A^{\ell,d}}} \right) \times 100 \%
\end{equation}
where $d$ is the test day, $O^{*,\ell,d}$ is the the profile after firming $A^{\ell,d}$ and for any $Z^{\ell,d}=(Z^{\ell,d}_{k})$,
$
\Dev{Z^{\ell,d}} := \sum_{k=0}^{23} \left| Z^{\ell,d}_k- F^{\ell,d}_k \right|
$
is the total daily deviation of the profile $Z$ from target generation ratio.
We proceed by training our algorithm for 24 randomly selected test days, 2 days from each month. Since we have $149$ hybrid assets, this totals to $149 \times 24$ runs of SHADOw-GP, each run taking $\approx$ 5 minutes. Left panel of Figure \ref{fig:L1reduction} visualizes the averaged $\overline{\mathrm{DR}}^{\ell} = \frac{1}{24} \sum_{d=1}^{24} \mathrm{DR}^{\ell,d}$ over the test days for the $149$ wind assets in Texas-7k. We observe $\overline{\mathrm{DR}}^{\ell}$ as low as $40\%$ and as high as $67\%$, with a clear spatial pattern whereby higher reductions are observed in the southern and southeastern Texas, as indicated by the brighter yellow bubbles, which also correspond to larger hybrid assets. We observe lowest $\overline{\mathrm{DR}}^{\ell}$  (\%) in the Far North  region of ERCOT. 

\subsection{Benefits of Hybridization}\label{ssec:4.3}
\begin{figure}[!ht]

    \centering
    \includegraphics[width = 0.8\linewidth]{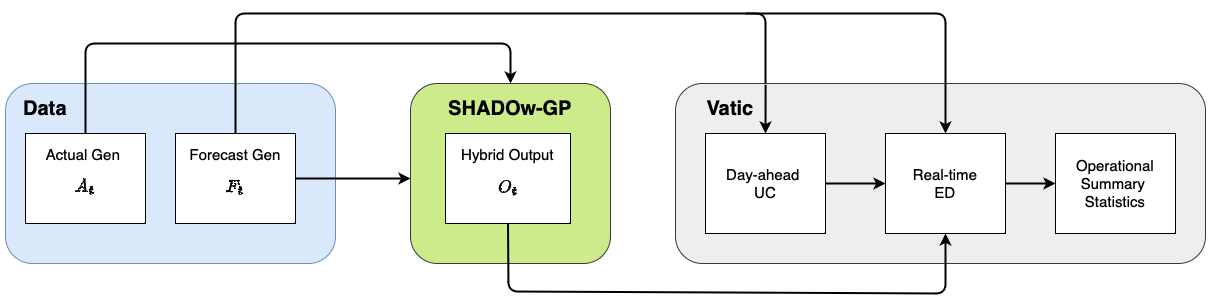}
    \caption{Integration of {SHADOW-GP} for Economic Dispatch in \texttt{Vatic}. }
    \label{fig:shadow-gp}
\end{figure}

To assess the economic benefits of retrofitting wind assets with BESS, we have  integrated the SHADOw-GP algorithm with the grid simulation tool \texttt{Vatic} \parencite{Vatic}. \texttt{Vatic} models daily grid operations via a two-step process depicted in Figure~\ref{fig:shadow-gp}: a day-ahead Unit Commitment (UC) step that schedules generators based on forecasts $F_k^{\ell,d}$ for $\ell=1,\ldots,149$, followed by an Economic Dispatch (ED) step that allocates generation in real time to minimize costs (under a 20\% reserve margin). For each day, SHADOw-GP takes the forecast and actual generation ratios as inputs and outputs firmed generation $O^{*,\ell,d}_k$, which \texttt{Vatic} then uses to run ED and compute operational cost savings.

For each hybrid asset $\ell$, we run \texttt{Vatic} with hybrid output $(O^{*,\ell,d}_k)$ for the same $D=24$ test days in the previous section, resulting in $149\times 24$ individual runs. Here, we emphasize that we turn only one wind asset into a hybrid for each run, with the same BESS characteristics as in Section \ref{ssec:5.2}. 
The right panel of Figure \ref{fig:L1reduction} reports the resulting ED savings in dollars per MW for each asset $\ell$, 
\begin{equation}
\text{ED savings}^{\ell} := \frac{ \sum_{d=1}^{D} \left( \text{ED}^{\text{Base}}_{d} - \text{ED}^{\text{Hybrid}, \ell}_{d} \right) }{ D \cdot B_{\max} \times \bar{C}_\ell }
\end{equation}
\noindent 
where
$\text{ED}^{\text{Base}}_{d}$ is the daily variable generation cost (\$) of the grid in the base run without hybrid assets (using solely $(A^{\ell,d}_k)$ rather than any $(O^{*,\ell,d}_k)$) on day $d$ and
$\text{ED}^{\text{Hybrid},\ell}_{d}$ is the cost for the grid with asset $\ell$ hybridized. To account for asset nameplate capacities, we normalize by 
$B_{\max} \times \bar{C}_\ell$ the power rating (MW) of the battery in asset $\ell$.

We observe average savings of approximately \$30/MW across the hybrid assets. While a few sites display values near \$140/MW, these correspond to isolated outliers rather than a systematic pattern. Supplementary Figure~S6 shows that larger savings are not consistently associated with small hybrid assets; rather, the spatial distribution reflects a mix of factors including grid location and local economic dispatch conditions. 


\subsection{BESS duration}
BESS duration is among the key parameters affecting firming performance. Larger energy capacity allows the hybrid asset to sustain target generation level during a longer deviation. Typical durations, which in our model correspond to the ratio $I_{\max}/B_{\max}$, are 1--6 hours, with 2-hour and 4-hour being the most common configurations for Lithium-ion battery storage circa 2025. 
To quantify the benefit to the grid of installing a longer-duration battery we compare  firming performance when we change from the $3$-hour BESS ($\Imax/B_{\max} = 3 $) of the previous section to a larger $6$-hour BESS ($\Imax/B_{\max} = 6$). 

The right panel of Figure~\ref{fig:L1reduction} visualizes the gains in terms of higher DR (more firming) from using a longer duration battery. We re-run \texttt{Vatic} ED, making $149 \times 24$ new runs as each potential hybrid retrofit is individually adjusted to $I_{\max}=0.6$. We observe gains of 6-13\% with longer-duration batteries. The intuition would be that gains are linked to the frequency of situations where more energy capacity is needed to charge/discharge over 3+ hours due to the particular  wind  generation remaining consistently above/below the forecast for several successive hours.


\section{Alternative objectives}\label{sec:6}

\subsection{Greedy Objective}

The quadratic firming criterion \eqref{eq:tildeF} used so far 
leads to a smooth relationship between $B_k(X, I)$ and SoC $I$. In contrast, a greedy controller reacts myopically by charging or discharging immediately, without any foresight. 
To illustrate, consider the $L_1$ objective:
\begin{equation}\label{eq:fgreedy}
    f_{L_1}(X_k,B_k,M_k) := |X_k-B_k-M_k|.
\end{equation}

In the absence of State-of-Charge constraints, the optimal action minimizes the instantaneous dispatch deviation, with no look-ahead \parencite{exarchos_stochastic_2018,nagahara2013L1optimality}. 
Formally, the sub-gradient optimality condition of the Bellman equation under \eqref{eq:fgreedy} is
\begin{equation}\label{eq:subgradientFOC}
    0 \in - \operatorname{sgn}(X_k-B_k-M_k) + 
    \frac{\partial \widehat{Q}_k}{\partial B_k} \!\left(X_k,\, I_{k}+
    \Bigl(\eta B_{k} \mathbb{1}_{\{B_k>0\}}+\tfrac{1}{\eta}B_{k}\mathbb{1}_{\{B_k<0\}}\Bigr)\Delta t\right).
\end{equation}

\begin{figure}[!ht]
    \centering
    \includegraphics[width=0.485\linewidth]
    {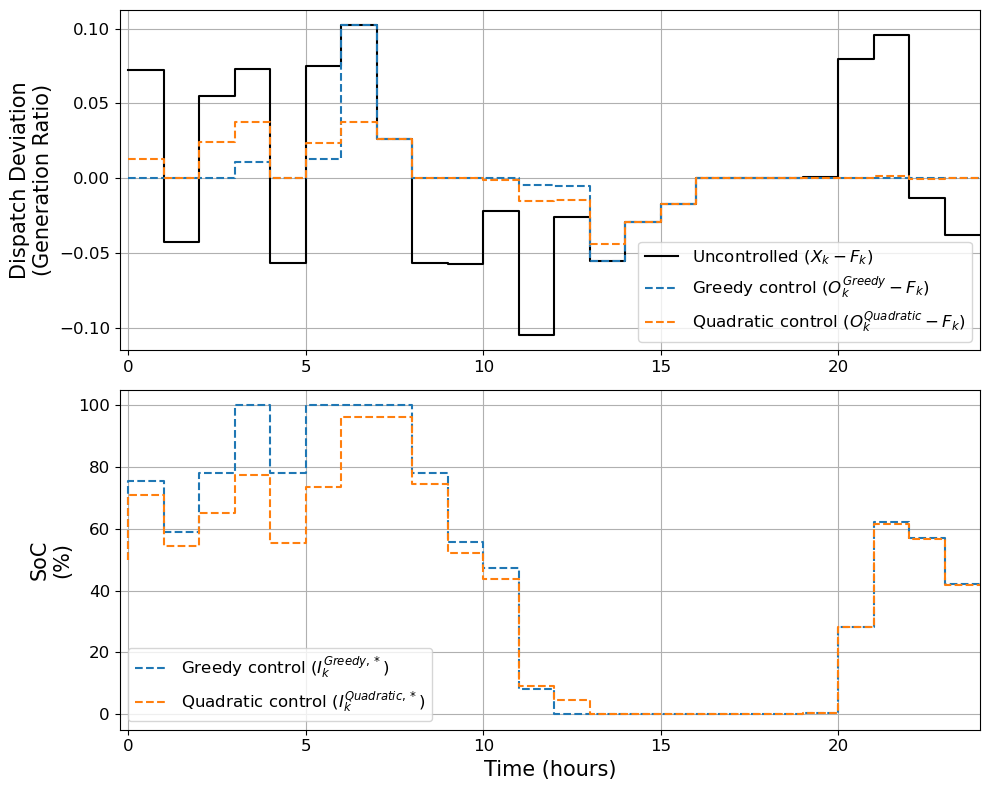}
    \includegraphics[width=0.485\linewidth]{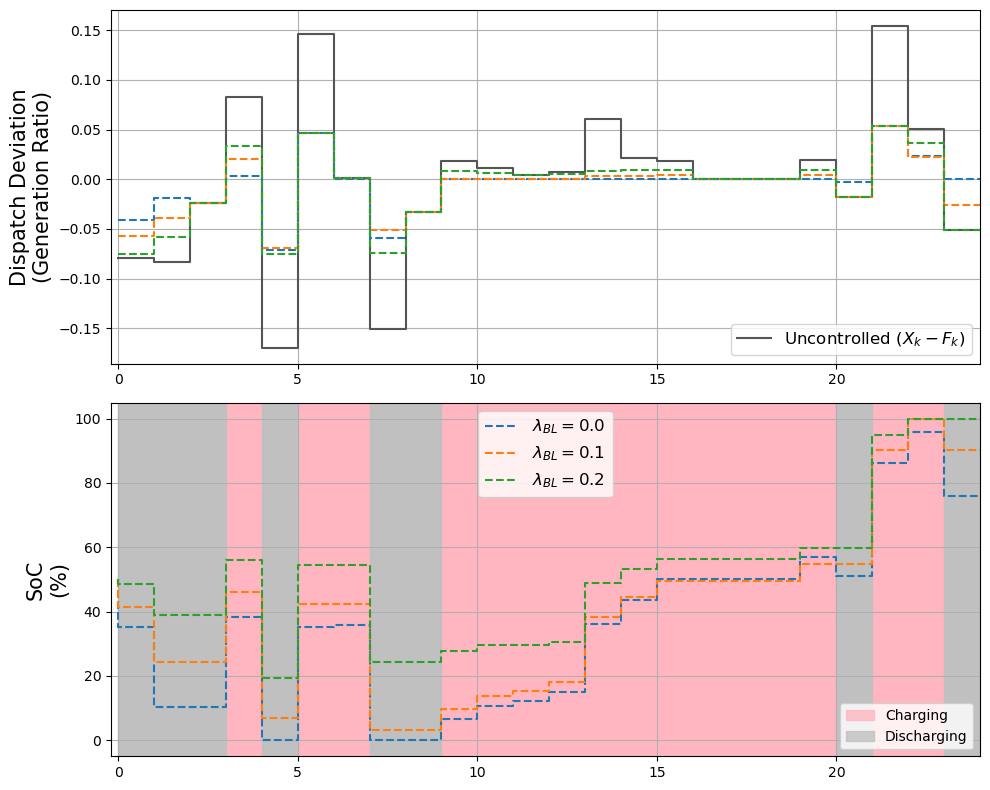}
    \caption{\emph{Top left panel:} A representative trajectory of dispatch deviation without firming $(X_k-F_k)$ and with firming $(O_k-F_k)$, the latter based on $\tilde{f}$ (quadratic control) or $f_{L_1}$ (greedy control). \emph{Bottom left panel}: Corresponding SoC trajectories $(I^{*}_k)$ optimizing $\tilde{f}$ and $f_{L_1}$.  \emph{Top  right panel:} A representative trajectory of dispatch deviation without firming $(X_k-F_k)$ and with firming $(O_k-F_k)$ over one day $k=0,\ldots, 23$. \emph{Bottom right panel:} Corresponding firmed SoC trajectories $I^{*}_k$. We show three scenarios based on battery lifetime penalty $\lambda_{BL} \in \{0, 0.1, 0.2\}$.
 }\label{fig:Greedy}
\end{figure}

This condition leads to three regimes:
\begin{itemize}
    \item If $\Big|\frac{\partial \widehat{Q}_k}{\partial B_k}(\cdot)\Big| <= 1$, 
    a stationary point exists at $B_k = X_k - M_k$. 
    In this case, the battery exactly cancels the deviation, subject to state and control limits.
    
    \item If $\frac{\partial \widehat{Q}_k}{\partial B_k}(\cdot) > 1$ (future cost of storage is high), 
    it is optimal to maximally discharge: $B_k = B_{\min}(I_k)$.
    
    \item If $\frac{\partial \widehat{Q}_k}{\partial B_k}(\cdot) < -1$ (future cost of storage is low), 
    the optimizer chooses full charge: $B_k = B_{\max}(I_k)$.
\end{itemize}
In our results, the first regime dominates: the optimal control cancels deviations exactly, confirming the greedy, no-look-ahead structure. To demonstrate this numerically, we run the SHADOw-GP algorithm for cost criterion  \eqref{eq:fgreedy}, for the case study of Section \ref{ssec:5.2}, utilizing the same BESS parameters.  The non-smooth absolute-value in \eqref{eq:fgreedy} causes the objective function in \eqref{eq:optimal_unconstrained} to be nonsmooth, increasing the run time by about a factor of $\simeq 4$ compared to the original quadratic objective, as the optimizer must numerically approximate the gradient. The left panels of Figure \ref{fig:Greedy} compare the paths of the re-centered firmed output and the SoC arising from the $L_1$ \eqref{eq:fgreedy} and quadratic \eqref{eq:tildeF} objectives. Evidently, the $L_1$ objective incites the BESS to swiftly charge or discharge to promptly minimize deviations, in a \emph{greedy} (myopic) fashion. 
Moreover, ${L_1}$-costs cause the BESS to frequently be 0\% or 100\% full, restricting flexibility at later steps, while the quadratic objective tries to strategically maintain headroom. 

\subsection{Degradation-aware Firming}
BESS cycling accelerates capacity fade and efficiency loss, reducing both usable energy capacity $\mathrm{SoC}_{\max}$ and roundtrip efficiency $\eta$ over time. To capture these effects in a tractable way, we follow \textcite{Heyman_degradation} and introduce a nonlinear degradation cost that proxies depth-of-discharge (DoD) through the state of charge and (dis)charge rates given by:
\begin{align}\label{eq:Phi}
\Phi (B,I)  =  \left(1- 0.5 (\frac{I}{I_{\max}})^2 \right)\cdot  \max(-B,0),
\end{align}
leading to the modified stepwise cost function: 
\begin{equation}\label{eq:fbl}
    f_{BL}(X_k,B_k,M_k) :=\tilde{f}(X_k,B_k,M_k) + \lambda_{BL} \cdot \Phi (B_k,I_k),
\end{equation}
where $\lambda_{BL}$ is the weight of the Battery Lifetime penalty. The second term in \eqref{eq:Phi} increases as $I \downarrow I_{\min}$ or $B \downarrow B_{\min}$ penalizing discharges of nearly empty battery that are associated with high DoD and promoting more shallow cycling.

To evaluate the effect of our degradation proxy $\Phi$, we compute degradation ex post over the optimized SoC trajectory $(I_k^*)$ obtained under $f_{BL}(\cdot)$. The SoC trajectory is decomposed into charging and discharging half-cycles via rainflow counting \parencite{yuanyuan}, and the daily degradation for asset $\ell$ on day $d$ is
\begin{equation}
L^{\ell,d}(\mathbf{c}, \mathbf{d}) = \sum_{i=1}^{|\mathbf{c}|} \frac{\Theta(c_i)}{2} + \sum_{i=1}^{|\mathbf{d}|} \frac{\Theta(d_i)}{2},
\end{equation}
where $\mathbf{c}$ and $\mathbf{d}$ are the depth-of-charge and depth-of-discharge vectors corresponding to $I_k^*$. The cycle degradation function we use is
\begin{equation}\label{eq:theta}
\Theta(\text{DoD}) = 5.24 \times 10^{-4}\,(\text{DoD})^{2.03},
\end{equation}
implying a lifetime of roughly $5.2$ years for daily full cycles.

The right panels of Figure~\ref{fig:Greedy} illustrate a representative trajectory behavior under $f_{BL}$ for $\lambda_{BL} \in \{0, 0.1, 0.2\}$ in the case study of Section~\ref{ssec:5.2}. As expected, larger $\lambda_{BL}$'s dampen charging and discharging responses to dispatch deviations, leading to smoother SoC trajectories and larger $(O_k - F_k)$ deviations.  

To assess this trade-off systematically, we deploy SHADOw-GP over $\lambda_{BL} \in \{0,0.02,0.04,\ldots,0.2\}$. We define the \emph{Expected Battery Life} (EBL) and \emph{Expected Deviation Reduction} (EDR) metrics as
\begin{equation}
\mathrm{EBL}^{\ell,d} := \EE\big[L^{\ell,d}(\mathbf{c},\mathbf{d})^{-1}\big],
\qquad
\mathrm{EDR}^{\ell,d} := \EE \left[ \frac{\Dev{X^{\ell,d}} - \Dev{O^{\ell,d}}}{\Dev{X^{\ell,d}}} \right] \times 100\%,
\end{equation}
capturing battery lifetime extension and firming performance, respectively.
The left panel of Figure~\ref{fig:degradation_tradeoff} shows a nearly linear dependence of $\mathrm{EBL}$ and $\mathrm{EDR}$ on $\lambda_{BL}$, computed from $10^4$ Monte Carlo simulations for Foard City Wind Farm on 2018-04-05. Increasing $\lambda_{BL}$ from $0$ to $0.2$ more than doubles the expected BESS lifetime (from $\sim$5 to $\sim$12 years) due to shallower cycles, at the cost of about a 25\% reduction in firming performance. This linear trend indicates that the proxy cost $\Phi$ closely mirrors rainflow-based degradation $L^{\ell,d}(\cdot)$ without requiring additional state variables for DoD tracking.

\subsection{Curtailment Mitigation}

As a final modification, we consider a controller that jointly minimizes the firming objective \eqref{eq:tildeF} and the asymmetric \emph{Expected Cumulative Violation} (ECV) relative to curtailment thresholds $(\overline{M}^{\ell,d}_k)$,
\begin{equation}\label{eq:tvc}
\mathrm{ECV}^{\ell,d} = \EE\left[\sum_{k=0}^{23} \max\bigl(X^{\ell,d}_k - B^{\ell,d}_k - \overline{M}^{\ell,d}_k, 0\bigr)\right],
\end{equation}
which proxies the aggregate curtailed energy over the day. The controller directs generation above $\overline{M}_k$ into the BESS while continuing to firm toward $M_k \le \overline{M}_k$, leading to the blended \emph{Curtailment Mitigation} (CM) stepwise cost with weight $\lambda_{CM}$
\begin{equation}
f_{\text{CM}}(X_k,B_k,M_k) := \tilde{f}(X_k,B_k,M_k) 
+ \lambda_{CM}\,\max\bigl(X_k - B_k - \overline{M}_k,0\bigr).
\end{equation}
Curtailment mitigation is critical in renewable-rich grids, where excess generation and limited downward regulation often lead to forced output reductions, wasted energy, and added costs; see \textcite{curtailment_GA} for a dynamic programming approach to wind curtailment. The thresholds $(\overline{M}^{\ell,d}_k)$ may also be constant, i.e., $\overline{M}^{\ell,d}_k \equiv \overline{M}^{\ell,d}$, representing fixed POI capacity limits. Many hybrid assets are oversized relative to interconnection capacity \parencite{LBNL}, so that $\overline{M}^{\ell,d}$ is below the nameplate rating and the BESS must absorb excess generation when outputs approach this limit.

\begin{figure}[!ht]
    \centering
    \includegraphics[width=0.475\linewidth]{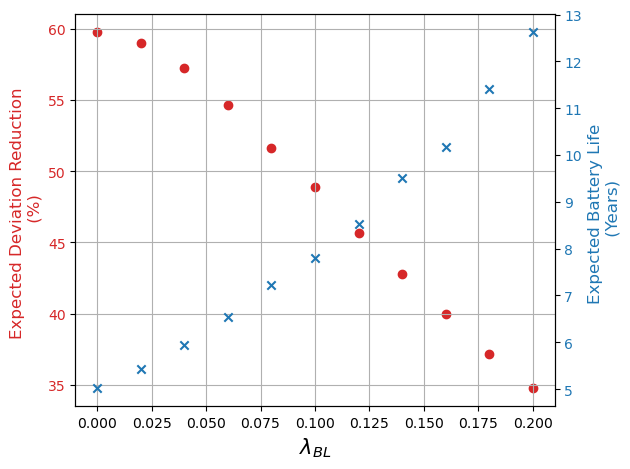}
    \includegraphics[width=0.475\linewidth]{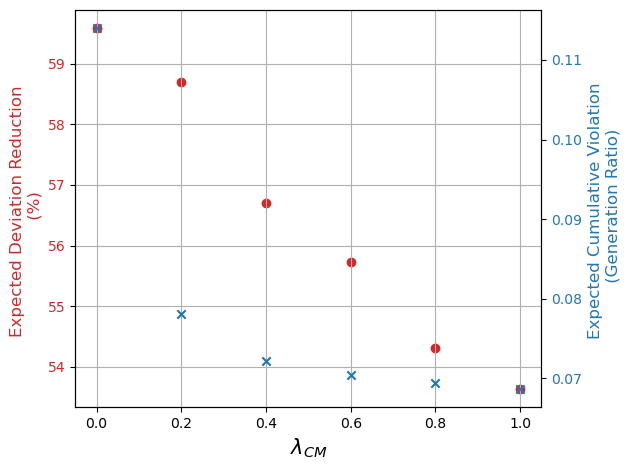}
    \caption{ \emph{Left panel:} Trade-off between $\mathrm{EDR}$ and $\mathrm{EBL}$, evaluated for $\lambda_{BL} \in \{0,0.02,0.04,\ldots,0.2\}$ for objective BL in \eqref{eq:fbl}.
    \emph{Right panel:} Trade-off between $\mathrm{EDR}$ and $\mathrm{\mathrm{ECV}}$, evaluated for $\lambda_{CM} \in \{0,0.2,\ldots,1.0\}$ for objective CM in \eqref{eq:tvc}. For both panels, the metrics are computed via $10^4$ Monte Carlo trajectories on Foard City Wind Farm, 2018-04-05.
   }
    \label{fig:degradation_tradeoff}
\end{figure}

The right panel of Figure~\ref{fig:degradation_tradeoff} illustrates the trade-off between firming minimization and curtailment cap exceedance, $\mathrm{EDR}$ and $\mathrm{ECV}$, across a range of weights $\lambda_{CM} \in [0,1]$, with cap set at $\overline{M}_k := 1.05\,M_k$ (5\% above the dispatch target). The $\mathrm{ECV}$ response to $\lambda_{CM}$ is highly non-linear, with a sharp decline between $\lambda_{CM}=0$ and $0.2$. Setting $\lambda_{CM}=1$ reduces $\mathrm{ECV}$ by about $40\%$ relative to $\lambda_{CM}=0$, while incurring only a $6\%$ loss in $\mathrm{EDR}$, indicating that the $f_{\text{CM}}$ criterion effectively mitigates curtailment with minimal firming performance decrease.

\section{Conclusion and Future Work}\label{sec:8}
 This paper presents a comprehensive and innovative framework for the optimal intraday dispatch of hybrid wind assets under multiple performance objectives, including quadratic penalty firming, degradation-aware operation, and curtailment mitigation. We develop a closed-loop algorithm for dynamic real-time dispatch that operates in continuous state and action spaces, is agnostic to the underlying wind dynamics $(X_k)$, and can leverage nonparametric scenario simulators. The accuracy of the solution is validated against an analytical benchmark in the linear–quadratic case. From an experimental perspective, we quantitatively assess the benefits of retrofit hybridization using a realistic Texas-7k grid synthetic testbed.

Our proposed algorithm controls the BESS under the assumption of a constant penalty for deviation from $(M_k)$. In practice, the penalty for missing the dispatch targets $(M_k)$ depends on real-time electricity prices $(P_k)$. Therefore, a hybrid controller might employ a time-sensitive firming strategy where the firming goal aligns with $P_k$. A natural future extension would be to undertake price-dependent firming by introducing an additional state variable $(P_k)$ to \eqref{eq:Bellman_DP}:
\begin{equation} 
             V(t_k,X_k,I_k,P_k ) = \inf_{B_k\in \mathcal{B}_k}  
\Bigl\{f(X_k,B_k,M_k,P_k) + \EE \left[V(t_{k+ 1} , X_{k+1},I_{k+1}(B_k), P_{k+1})|{X_k,I_k,P_k} \right] \Bigr\}
\end{equation}
where the criterion $f$ depends on $p$, such as $f(x,b,m,p) = p(x-b-m)^2$.

The presented model focuses solely on real-time dispatch tracking of hybrid assets. In a different vein, assets can be active in both the day-ahead (DA) and real-time (RT) markets. Asset dispatchers can first strategically participate in the DA market, and then dynamically adjust their positions via the RT market intraday. This leads to a two-step optimization in which RT adjustments are impacted by earlier DA commitments, and can be recast as a firming objective where the target $(M_t)$ is itself optimized via a further objective based on the DA prices, see e.g., \textcite{tankov_wind_BESS}. Full co-optimization  requires capturing the probabilistic dependence between DA and RT prices.

From a control perspective, an important generalization would be to adapt our control emulator to a multivariate setting with constraints, which would open the door to study BESS use-case stacking, such as the hybrid asset jointly participating in firming and ancillary service provision.


\section{Acknowledgements}
Work partially supported by NSF-2420988 and NSF-2407550. We thank Haosheng Zhou and Samuel Babichenko for useful discussions. 
\printbibliography

\clearpage

\section*{Supplementary Materials}

\beginsupplement
\section{Proofs }\label{app:hjb}
\noindent \textbf{Proposition 1.} \textit{Proof:} 
\noindent The unconstrained HJB equation corresponding to the linear–quadratic control problem presented in the main paper is given by:
\begin{align}\label{HJB_proof}
\partial_t \bar{V} 
+ \inf_{\bar{b}\in \mathcal{R}} \Bigl\{
& \alpha_t (m - x)\,\partial_X \bar{V} 
+ \bar{b}\,\partial_I \bar{V} 
+ \tfrac{1}{2}\sigma_t^2\, x\,(\Xmax - x)\,\partial_{XX} \bar{V} \nonumber\\
& \quad + (x - \bar{b} - M_t)^2 
+ c_1\,\bar{b}^2 
+ c_2\,(\iota - I_m)^2
\Bigr\} = 0.
\end{align}

with terminal condition  $\bar{V}(T,x,\iota) = \mathcal{P} (\iota - I_{\text{target}})^2$. Since the value function is quadratic in $x, \iota$ with interaction terms, we can derive an analytical formulation for the optimal control and $\bar{V}$. 

The first-order condition of $\bar{b}$ in \eqref{HJB_proof} is:
\begin{equation}\label{bbar}
    \bar{b} = \kappa (x-M_t) - \frac{\kappa}{2} (\partial_{I} \bar{V}),
\end{equation}
where $\kappa = \frac{1}{1+c_1}$.
Plugging back into \eqref{HJB_proof} yields
\begin{align*}
    &\partial_t \bar{V} 
    + \alpha_t(m-x) \partial_X \bar{V} 
    + \left[\kappa (x-M_t) - \frac{\kappa}{2} \partial_{I} \bar{V}\right] \partial_{I} \bar{V} 
    + \frac{1}{2}\sigma_t^2 x \cdot(\Xmax - x)\partial_{XX} \bar{V} 
     + \\ & \quad \left[(x-M_t) - \kappa (x-M_t) + \frac{\kappa}{2} \partial_{I} \bar{V} \right]^2 
     + c_1 \left[\kappa (x-M_t) - \frac{\kappa}{2} (\partial_{I} \bar{V})\right]^2 
    + c_2 (\iota-I_m)^2 = 0.
\end{align*}
Expanding and simplifying  we end up with 
\begin{align*}
\partial_t \bar{V} 
&+ (1 - \kappa)\,(x - M_t)^2 
+ c_2\,(\iota - I_m)^2 
- \alpha_t\,(x - m)\,\partial_X \bar{V} \nonumber\\
&\quad + \kappa\,(x - M_t)\,\partial_I \bar{V} 
- \frac{\kappa}{4}\,\bigl(\partial_I \bar{V}\bigr)^2 
+ \frac{1}{2}\,\sigma_t^2\,\partial_{XX} \bar{V} 
= 0.
\end{align*}
We rewrite $x-M_t = (x-m) + (m-M_t)$ to yield:
\begin{align*}
\partial_t \bar{V} &+ (1 - \kappa)(x - m)^2 + 2(1 - \kappa)(m - M_t)(x - m) + (1 - \kappa)(m - M_t)^2 \\
&+ c_2(\iota - I_m)^2 - \alpha_t(x - m)\partial_X V + \kappa(x - m)\partial_I \bar{V} + \kappa(m - M_t)\partial_I V \\
&- \frac{\kappa}{4}(\partial_I \bar{V})^2 + \frac{1}{2}\sigma_t^2 x \cdot(\Xmax - x)\partial_{XX} \bar{V} = 0.
\end{align*}

We solve the PDE by proposing the ansatz
\begin{equation}\label{v_ansatz}
\begin{aligned}
\bar{V}(t, x, \iota) &= P_1(t)(\iota - I_m)^2 + P_2(t)(\iota - I_m)(x - m) + P_3(t)(x - m)^2 \\
&\quad + P_4(t)(\iota - I_m) + P_5(t)(x - m) + P_6(t)
\end{aligned}
\end{equation}
for functions of time $P_j(t), j=1,\ldots, 6$ to be determined. Substituting back and gathering terms we get

\begin{align*}
&\dot{P_1}(t)(\iota - I_m)^2 + \dot{P_2}(t)(\iota - I_m)(x-m) + \dot{P_3}(t)(x-m)^2 + \dot{P_4}(t)(\iota - I_m) + \dot{P_5}(t)(x-m) + \dot{P_6}(t) \\
&+ (1-\kappa)(x-m)^2 + 2(1-\kappa)(m-M_t)(x-m) + (1-\kappa)(m-M_t)^2 + c_2(\iota - I_m)^2 \\
& - \alpha_t(x-m)\left[P_2(t)(\iota - I_m) + 2P_3(t)(x-m) + P_5(t)\right] \\
& + \kappa(x-m)\left[2P_1(t)(\iota - I_m) + P_2(t)(x-m) + P_4(t)\right] \\
& + \kappa(m-M_t)\left[2P_1(t)(\iota - I_m) + P_2(t)(x-m) + P_4(t)\right] \\
& - \frac{\kappa}{4}\left(2P_1(t)(\iota - I_m) + P_2(t)(x-m) + P_4(t)\right)^2   + \sigma_t^2P_3(t) x \cdot(\Xmax - x) = 0.
\end{align*}

Since the basis $\{(x-m)^2,(x-m)(\iota - I_m),(\iota - I_m)^2, (x-m),(\iota - I_m),1\}$ is linearly independent, we must have the respective coefficients identically zero, which yields the system of Riccati ODEs for $P_1(t),\ldots, P_6(t)$ on $t \in [0,T]$ below.

{\small
\begin{align} \left\{ 
\begin{aligned}
&\dot{P_1}(t) + c_2 - \kappa (P_1(t))^2 = 0, & P_1(T) = \mathcal{P}\\
&\dot{P_2}(t) -\alpha_t P_2(t) -\kappa P_1(t) P_2(t) + 2 \kappa P_1(t) = 0, & P_2(T) = 0 \\
&\dot{P_3}(t)+(1-\kappa) -(2 \alpha_t + \sigma^2_t) P_3(t) +\kappa P_2(t) -\frac{\kappa}{4} (P_2(t))^2 = 0, & P_3(T)=0 \\
&\dot{P_4}(t) + 2\kappa(m-M_t) P_1(t) - \kappa P_1(t) P_4(t) = 0, &  P_4(T) = 2\mathcal{P}(I_m - \iota_{\text{target}})\\
&\dot{P_5}(t) + 2(1-\kappa) (m-M_t)-\alpha_t P_5(t) + \kappa P_4(t) + \kappa(m-M_t) P_2(t)  \\
&\quad -\frac{\kappa}{2} P_2(t)P_4(t) + (\sigma^2_t \Xmax - 2 \sigma^2_t m) P_3(t)= 0, & P_5(T) = 0 \\
&\dot{P_6}(t) + (1-\kappa) (m-M_t)^2 + \kappa (m-M_t) P_4(t) -\frac{\kappa}{4} (P_4(t))^2 = 0, & P_6(T) = \mathcal{P}(I_m - \iota_{\text{target}})^2
\end{aligned} \right. 
\end{align}
}





\begin{remark} When $I_{target} = I_{m}$ and $\tilde{m}_k = m$, the linear terms in $\bar{V}$ disappear as $P_4$ and $P_5$ become degenerate. When $M_t = m$ only, the linear term with respect to $(x-m)$ in $\bar{V}$ disappears. However, in the case $\iota_{target} = \iota_{m}$ only, none of the terms simplify due to the presence of $\kappa(m - M_t)\partial_I V$ in the PDE.
\end{remark}

\noindent \textbf{Proposition 2.} \textit{Proof:} 

\noindent \textbf{Local existence and uniqueness.} 
Define 
\[
\mathbf{P}(t) = \big(P_1(t), P_2(t), P_3(t), P_4(t), P_5(t), P_6(t)\big)^{\top}.
\] 
Then the nonlinear Riccati system \eqref{eq:ricatti} can be written in compact form as  
\[
\dot{\mathbf{P}}(t) = \mathbf{F}(t,\mathbf{P}(t)),
\]
where $\mathbf{F}: D \subset \mathbb{R}\times \mathbb{R}^6 \to \mathbb{R}^6$, with $D$ a closed rectangle such that $(T,\mathbf{P}(T)) \in \operatorname{int} D$, and  
\[
\mathbf{P}(T) = \big(\mathcal{P}, 0, 0, 2\mathcal{P}(I_m-\iota_{\text{target}}), 0, \mathcal{P}(I_m-\iota_{\text{target}})^2\big)^T.
\]
The mapping $\mathbf{F}$ is continuous in $t$ and locally Lipschitz in $\mathbf{P}$ (since it is at most quadratic with cross–interaction terms in $P_i$). By the Picard–Lindelöf theorem, there exists $\epsilon > 0$ such that $\mathbf{P}(t)$ admits a unique solution on the interval $[T-\epsilon, T+\epsilon]$.

\noindent \textbf{Global existence and uniqueness.} We show that by solving $P_i(t)$ iteratively, starting from $P_1(t)$, the solution for $\mathbf{P}(t) $ exists and unique over $[0,T]$ for any $T>0$. Consider the scalar Riccati equation
\[
\dot P_1 = \kappa P_1^2 - c_2, 
\qquad P_1(T) = \mathcal P, 
\qquad \kappa,c_2 > 0.
\]
Let $\gamma := \sqrt{c_2/\kappa} > 0$ and set 
\[
r := \frac{\mathcal P - \gamma}{\mathcal P + \gamma}.
\]
A standard separation of variables yields the explicit solution
\[
P_1(t) \;=\; \gamma \,
\frac{1 + r\,e^{2\kappa \gamma (t-T)}}
     {1 - r\,e^{2\kappa \gamma (t-T)}}.
\]
Since $\mathcal P>0$ we have $\mathcal P+\gamma>0$, hence $r \in (-1,1)$. 
For all $t \in [0,T]$ one has $e^{2\kappa \gamma (t-T)} \in (0,1]$, so
\[
|r\,e^{2\kappa \gamma (t-T)}| < 1.
\]
Therefore the denominator never vanishes on $[0,T]$, and the solution 
$P_1(t)$ remains finite and continuous throughout the interval. 
Combined with local existence and uniquness from Picard-Lindelof Theorem, the existence and uniqueness of $P_1(t)$ can be extended to $[0,T]$ for any $T>0$.Now, we show global existence and uniqueness of $P_i, i =2,\ldots,6$. 

\noindent $P_2$ satisfies 
  \[
  \dot P_2 = \big(\alpha_t + \kappa P_1(t)\big) P_2 - 2\kappa P_1(t).
  \]
Since $P_1(t)$ is already available in closed form, the equation for $P_2$ is a first-order linear ODE with known coefficients. By the method of integrating factors it admits an explicit solution.

\noindent $P_3$ satisfies 
  \[
  \dot P_3 = \big(2\alpha_t + \sigma_t^2\big) P_3 - \kappa P_2(t) + \tfrac{\kappa}{4}P_2(t)^2 -(1-\kappa).
  \]
This is again first-order linear ODE in $P_3$ with $P_2$ already known. Similarly, $P_3(t)$ admits an explicit solution. 

\noindent $P_4$ satisfies 
  \[
  \dot P_4 = \kappa P_1(t) P_4 - 2\kappa\,(m-M_t)\,P_2(t).
  \]
  With $P_1,P_2$ already known, this is a linear ODE in $P_4$ with explicit integrating factor solution.
  
  \noindent $P_5$ satisfies 
  \[
  \dot P_5 = \alpha_t P_5 - \kappa P_4(t) - \kappa(m-M_t)P_2(t)
             + \tfrac{\kappa}{2}P_2(t)P_4(t) - \big(\sigma_t^2 X_{\max}-2\sigma_t^2 m\big) P_3(t)
             - 2(1-\kappa)(m-M_t).
  \]
  Given $P_2,P_3,P_4$, this is linear in $P_5$ and admits an explicit solution. Finally, $P_6$ satisfies 
  \[
  \dot P_6 = -\kappa(m-M_t)P_4(t) + \tfrac{\kappa}{4}P_4(t)^2 - (1-\kappa)(m-M_t)^2,
  \]
  which is an inhomogeneous linear ODE with no $P_6$ term. It integrates directly to
  \[
  P_6(t) = P_6(T) + \int_t^T \Big[-\kappa(m-M_u)P_4(u) + \tfrac{\kappa}{4}P_4(u)^2 - (1-\kappa)(m-M_u)^2\Big]\,du.
  \]

 Since $\alpha_t, \sigma_t, M_t \in L^1([0,T])$, the explicit solutions for $P_2(t),\dots,P_6(t)$ are well defined. Moreover, each $P_i(t)$ depends only on previously constructed components. Because $P_1(t)$ exists uniquely on $[0,T]$ for any $T>0$, it follows that $P_2(t),\dots,P_6(t)$ also admit unique solutions on $[0,T]$ for any $T>0$. Hence, the full Riccati system~\eqref{eq:ricatti} admits a unique global solution on $[0,T]$.
\section{ Texas-7k Case Study}\label{Appen_Texas7k}
\subsection{Topology of the Texas-7k Power Grid}

\begin{figure}[H]
    \centering
  \includegraphics[width=0.45\linewidth]{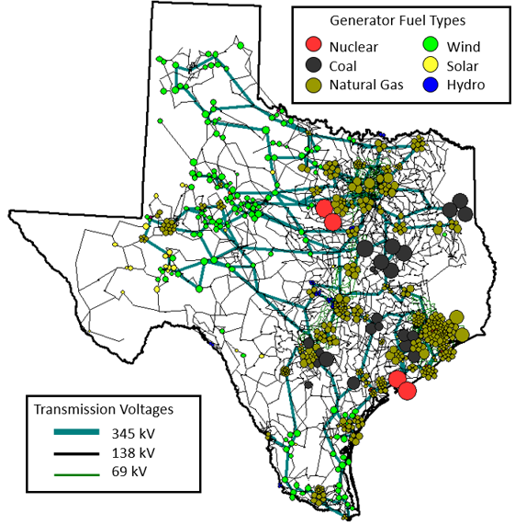}
    \caption{Topology of Texas-7k power grid. }
    \label{fig:texas7k_topology}
\end{figure}

Figure~\ref{fig:texas7k_topology} illustrates the topology of the Texas-7k test case, which represents a high-fidelity model of the ERCOT transmission network. The grid comprises 7{,}173 transmission lines and 1{,}967 transformers, with voltage levels spanning 345~kV, 138~kV, and 69~kV. The system includes a diverse mix of generation resources: 36 solar, 149 wind, and 542 thermal generators (coal, nuclear, natural gas, and hydro).

\subsection{Calibration of Wind Assets in Texas-7k}

\begin{figure}[H]
    \centering
    \includegraphics[width=0.65\textwidth]{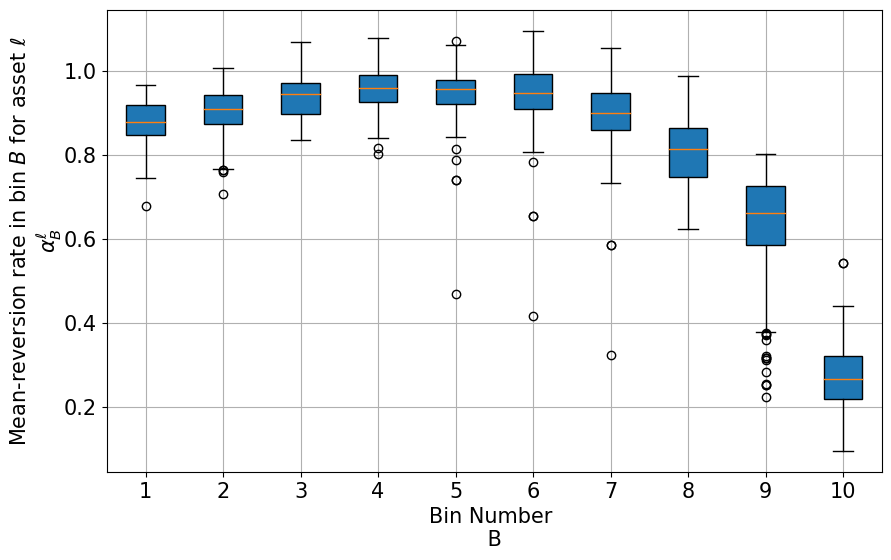}
    \caption{Boxplot of the calibrated mean-reversion rates $\alpha^{\ell}_{r}$ across 149 wind assets in Texas-7k as a function of bin $r=1,2,\ldots,10$. Mean-reversion is stable for moderate output levels and declines at high outputs.}
    \label{fig:alpha-b}
\end{figure}

In the main text, we introduced a bin-wise calibration of the drift and diffusion terms in the discrete-time mean-reverting dynamics of wind generation. Figure~\ref{fig:alpha-b} shows the boxplot of the calibrated mean-reversion rates $\alpha^{\ell}_{r}$ across 149 wind assets in the Texas-7k system, for bins $r=1,\ldots,10$. The results indicate that mean-reversion rates remain high and stable for moderate forecast levels (bins 1–7), but decline sharply in the upper bins. This reflects slower adjustment dynamics at high wind outputs, where persistent deviations from forecasts are common.

\begin{figure}[H]
    \centering
    \includegraphics[width=0.75\textwidth]{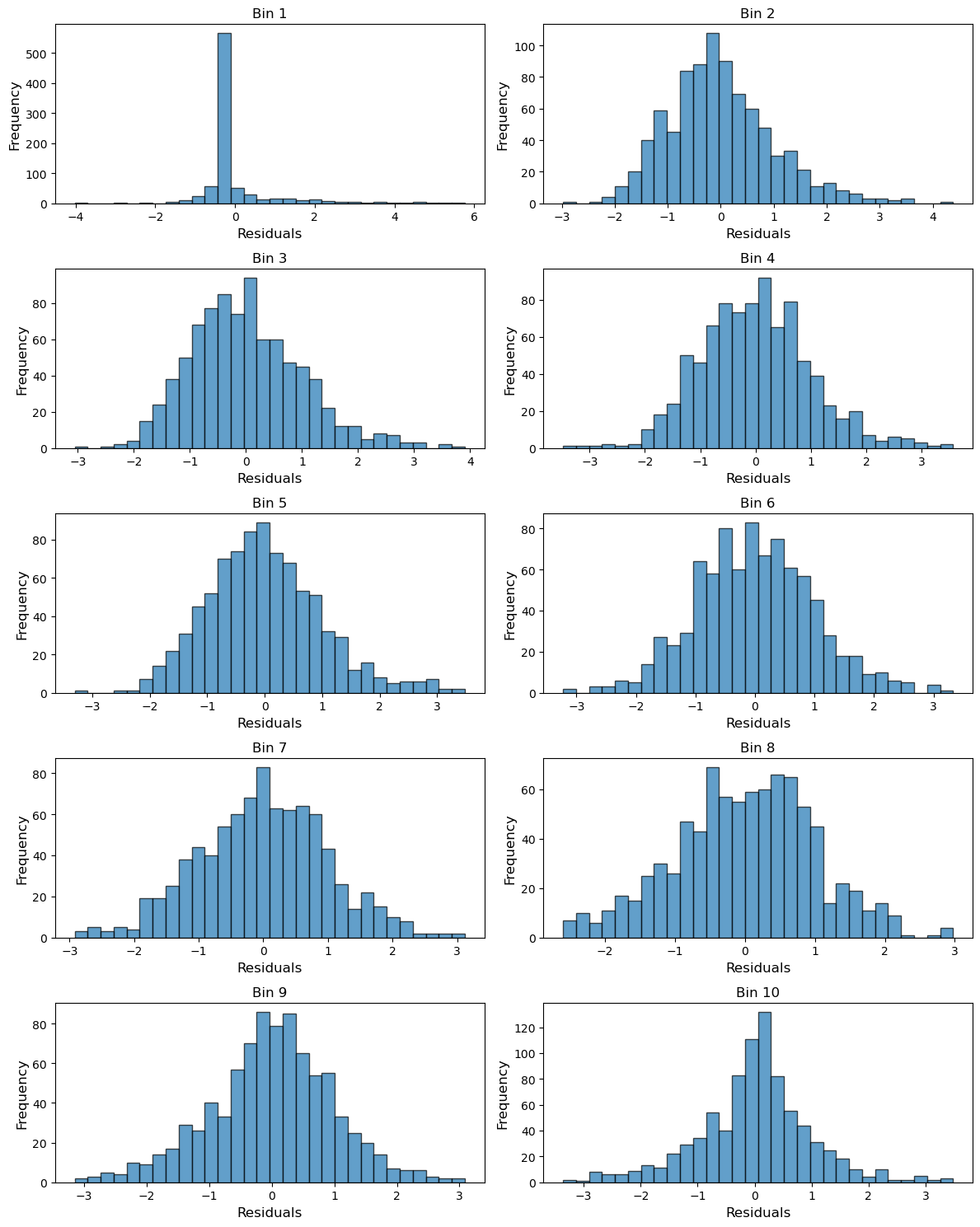}
    \caption{Histograms of bin-wise residuals $\mathcal{E}^\ell_r$ after calibrating $\alpha^\ell_r$ and $\sigma^\ell_r$ for the Foard City wind unit, across bins $r = 1,\dots,10$. The distributions exhibit skewness and heavy tails, deviating from Gaussianity.}
    \label{fig:residuals}
\end{figure}

Figure~\ref{fig:residuals} displays histograms of bin-wise residuals $\mathcal{E}^{\ell}_{r}$ after calibration for the Foard City wind unit. While residuals are centered around zero, their shapes deviate from Gaussianity in several bins, with heavier tails and skewness, particularly at intermediate and high-output regimes. This illustrates the presence of non-Gaussian features in the wind generation shocks and supports the use of more flexible noise models beyond simple Gaussian assumptions.

\newpage
\subsection{Empirical Coverage}
\begin{figure}[H]
\begin{center}
  \includegraphics[width=0.7\linewidth]{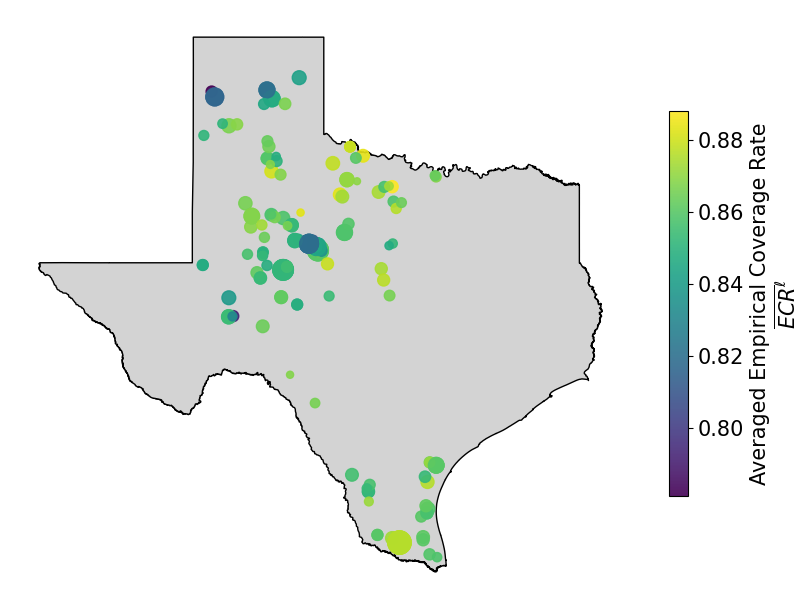}
    \caption{Spatial distribution of the averaged empirical coverage rate $\overline{\mathrm{ECR}}^{\ell}$ across 149 wind assets in the Texas-7k grid. Each point represents a wind asset, with symbol size proportional to its nameplate capacity and color indicating $\overline{\mathrm{ECR}}^{\ell}$.}
\label{fig:80CIsites}
\end{center}

\end{figure}

To check how well do the simulated trajectories $X^{\ell,d}_k$ match the observed actual generation ratios $A^{\ell,d}_k$ we compute the empirical time-averaged coverage rate $\mathrm{ECR}^{\ell,d}$ 
\begin{equation}\label{eq:ecr}
    \mathrm{ECR}^{\ell,d} := \frac{1}{24} \sum_{k=0}^{23} \mathbb{1} \Big( A^{\ell,d}_k\geq q_{0.1} (\mathbf{X}^{\ell,d}_k)\text{ and }   A^{\ell,d}_k \leq q_{0.9} (\mathbf{X}^{\ell,d}_{k}) \Big),
\end{equation}
where  $q_\alpha(\cdot)$ denotes the $\alpha$-th quantile. We take nominal coverage level of $\alpha=0.8$, corresponding to checking whether the 80\% band contains $A^{\ell,d}_k$. 
Figure \ref{fig:80CIsites} visualizes the averaged empirical coverage rate $\overline{\mathrm{ECR}}^{\ell}=\frac{1}{365} \sum_{d=1}^{365} \mathrm{ECR}^{\ell,d}$ across the Texas-7k assets. We find that $\overline{\mathrm{ECR}}$ ranges from 78.1\% to 88.8\%. This range is consistent with the nominal 80\% level, indicating overall good calibration with slight over-coverage for some assets and spatial variation reflecting heterogeneous wind characteristics.



\newpage
\subsection{Time-Dependent Structure of the Control Map}
\begin{figure}[H]
    \centering
    \includegraphics[width=0.6\linewidth]{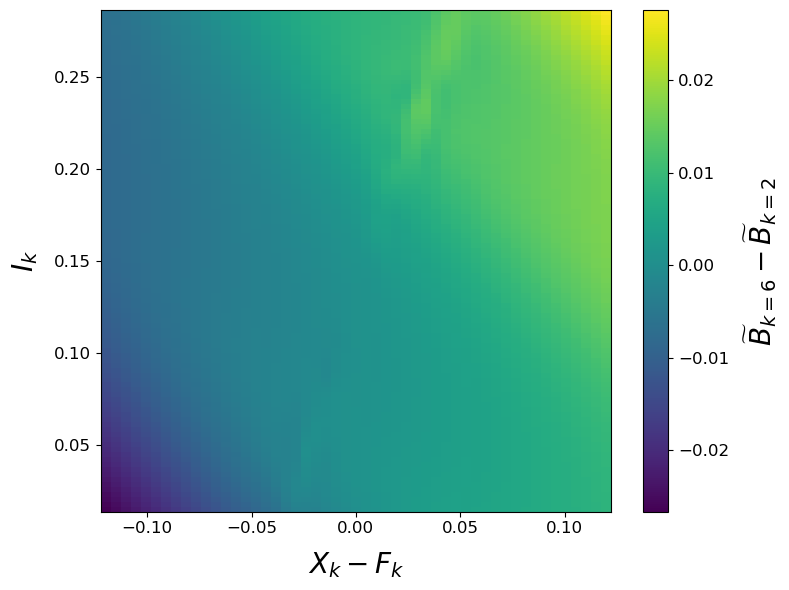} 
    \caption{Difference of unconstrained optimal control policies, $\tilde{B}_{k}(\cdot)$, from SHADOw-GP at different times, $k=6$ and $k=2$ as a function of $(X,I)$, calibrated to Foard City Wind Farm on 2018-04-05. For better visualization of deviation levels, we plot the dispatch deviation $X_k-F_k$ rather than the generation ratio $X_k$ itself. BESS charges/discharges more at hour $k=6$ than at $k=2$, demonstrating time-dependent control behavior.}
    \label{fig:stacked_control_map}
\end{figure}

The time-dependent nature of the underlying system dynamics naturally leads to time-dependent control policies. SHADOw-GP captures this non-stationarity by conditioning the optimal unconstrained control $\tilde{B}_{k}(X,I)$ on both the current state and the hour $k$. Figure~\ref{fig:stacked_control_map} illustrates this behavior by comparing control maps at two representative times.

\newpage
\subsection{Economic Dispatch (ED) of Hybrid Wind--Battery Assets}
\begin{figure}[H]
    \centering
    \includegraphics[width=0.7\linewidth]{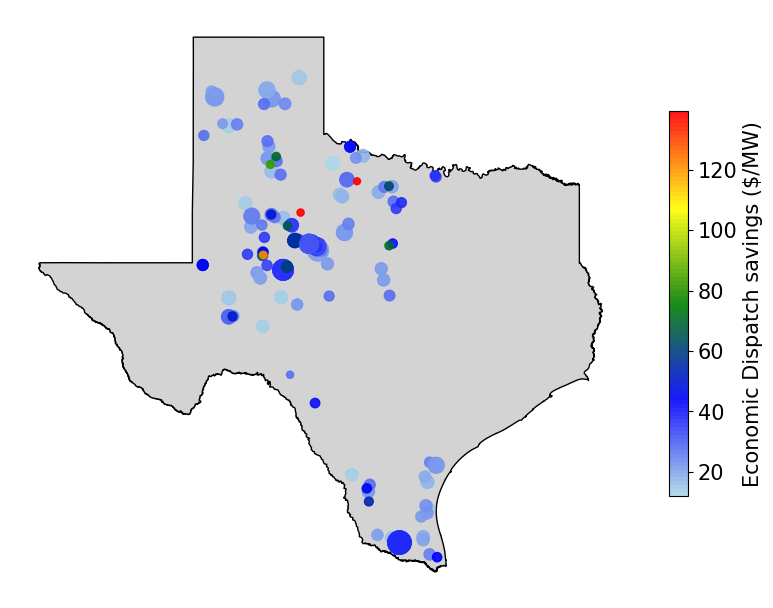}
    \caption{ED savings for wind asset retrofitted 3-hour BESS at 10\% capacity.}
    \label{fig:deviation_reduction_difference}
\end{figure}

Figure~\ref{fig:deviation_reduction_difference} summarizes the economic impact of retrofitting wind assets in the Texas-7k system with 3-hour BESS sized at 10\% of nameplate capacity. For each of the 149 wind assets, we evaluate average ED cost savings over 24 test days when hybrid dispatch is co-optimized using SHADOw-GP's firmed output. The color scale indicates savings in \$/MW, while marker size reflects asset capacity. 

\end{document}